\documentclass[twoside,english,10pt]{article}
\usepackage[latin9]{inputenc}
\usepackage{babel}
\usepackage[T1]{fontenc}
\usepackage[a4paper]{geometry}
\usepackage{amsthm}
\usepackage{amsmath}
\usepackage{amssymb,bm}
\usepackage{graphicx}
\usepackage{setspace}
\usepackage{esint}
\usepackage{epstopdf}
\usepackage{color}
\usepackage{appendix}

\def\Z{{\mathbb Z}}

\def\le{\leqslant}
\def\ge{\geqslant}

\newcommand{\im}{\mathrm{Im}}
\newcommand{\eps}{\varepsilon}
\newcommand{\e}{\varepsilon}

\begin{document}

\title{The Gaussian
wave packets transform  for the semi-classical Schr\"odinger
equation with vector potentials}

\author{Zhennan Zhou \thanks{Beijing International Center for Mathematical Research, Peking University, No. 5 Yiheyuan Road Haidian District,  100871, Beijing, China.}, \quad Giovanni Russo\thanks{Department of Mathematics and Computer Science, University of Catania, Via A.Doria 6, 95125, Catania, Italy.} }
\date{The paper is dedicated to the memory of Professor Peter Smereka}

\maketitle

\begin{abstract}
In this paper, we reformulate the semi-classical Schr\"odinger
equation in the presence of electromagnetic field by the Gaussian
wave packets transform.  With this approach, the highly oscillatory Schr\"odinger equation is equivalently transformed into another Schr\"odinger type wave equation, the $w$ equation, which is essentially not oscillatory and thus requires much less computational effort. We propose two numerical methods to solve the $w$ equation, where the Hamiltonian is either divided into the kinetic, the potential and the convection part, or into the kinetic and the potential-convection part. The convection, or the potential-convection part is treated by a semi-Lagrangian method, while the kinetic part is solved by the Fourier spectral method. The numerical methods are proved to be unconditionally stable, spectrally accurate in space and second order accurate in time, and in principle they can be extended to higher order schemes in time. Various one dimensional and multidimensional numerical tests are provided to justify the properties of the proposed methods.  
\end{abstract}

\section{Introduction}

In this paper, we propose an efficient approach for the semi-classical Schr\"odinger equation with external electromagnetic fields. This approach is a natural but worthy extension of the Gaussian Wave Packet Transform developed in \cite{GB Trans, GB Trans2} to the cases where the Hamiltonians include vector potentials. With this formulation, the highly oscillatory Schr\"odinger equation is transformed into another Schr\"odinger type wave equation, which is much smoother and, as a consequence, requires much less computational effort. This problem is challenging for various reasons, and has many important applications in physics and chemistry (see \cite{FL1,FL3,Lubichbook,reviewsemiclassical}). 

Consider the dimensionless Schr\"odinger equation for a charged particle, with a small (scaled)
Planck constant $\varepsilon$,
\begin{equation}
i\varepsilon\partial_{t}\psi^{\varepsilon}=\frac{1}{2}\left(-i\varepsilon\nabla_{\mathbf{x}}-\mathbf{A}(\mathbf{x})\right)^{2}\psi^{\varepsilon}+V(\bm{x})\psi^{\varepsilon},\quad t\in\mathbb{R}^{+},\quad\mathbf{x}\in\mathbb{R}^{3};\label{eq:main equation}
\end{equation}
\begin{equation}
\psi^{\varepsilon}(\mathbf{x},0)=\psi_{0}(\mathbf{x}),\quad\mathbf{x}\in\mathbb{R}^{3}\label{eq:initial cond}
\end{equation}
where $\psi^{\varepsilon}(\mathbf{x},t)$ is the complex-valued
wave function, $V(\mathbf{x})\in\mathbb{R}$  is the scalar potential and $\mathbf{A}(\mathbf{x})\in\mathbb{R}^{3}$
is the vector potential. The scalar potential  and the vector
potential  are introduced to mathematically describe  the external electromagnetic field, or respectively, the electric
field $\mathbf{E}(\mathbf{x})\in\mathbb{R}^{3}$ and the magnetic
field $\mathbf{B}(\mathbf{x})\in\mathbb{R}^{3}$ given by,
\begin{equation}
\mathbf{E}=-\nabla V(\mathbf{x}),\qquad\mathbf{B}=\nabla\times\mathbf{A}\left(\mathbf{x}\right).
\end{equation}
For simplicity, we have assumed that the scalar and the vector potentials are time independent, but inclusion of the time dependence will not add  intrinsic challenges to the current problem. 

The quantum Hamiltonian in \eqref{eq:main equation} is reminiscent of the classical Hamiltonian  (see \cite{Jacksonbook})  
\[
H(\bm q, \bm p) = \frac 1 2 (\bm p - \bm A(\bm q) )^2 + V(\bm q).
\]
This is the  Schr\"odinger equation for a charged quantum particle moving in an electromagnetic field, where no spin or relativistic effects are considered (see \cite{FL3}). It also shows a connection between quantum mechanics and macroscopic scale effects (the classical electromagnetic fields in this context). Alternatively, this equation (\ref{eq:main equation})
above can be derived from the free-particle Schr\"odinger  equation 
by local gauge transformation (see \cite{QO}). 


The quantum dynamics in the presence of  external
electromagnetic fields results in many far-reaching consequences in
quantum mechanics, such as Landau levels, Zeeman effect and superconductivity (see, e.g., \cite{FL3}).
Mathematically, it gives new challenges as well,
especially in the semi-classical regime. The presence of the vector
potential introduces a convection term in the Schr\"odinger
equation and in the meanwhile effectively modifies the scalar potential
(see \cite{SL-TS}). 

%

In the semi-classical regime, namely $\varepsilon \ll 1$, the solution to the Schr\"odinger
equation is highly oscillatory both in space and time on the scale
$O(\varepsilon)$, so that the wave function $\psi^{\varepsilon}$ does
not converge in the strong sense as $\varepsilon\rightarrow0$. Numerically,
the oscillatory nature of the wave function of the semi-classical
Schr$\rm{\ddot{o}}$dinger equation gives rise to significant computational
burdens. If one aims for direct simulation of the wave function, according to our knowledge, one of the best choices is the time splitting spectral method introduced by Bao, Jin and Markowich
in \cite{TSSP,reviewsemiclassical}, where the meshing strategy $\Delta t=O(\varepsilon)$
and $\Delta x=O(\varepsilon)$ is sufficient for moderate values of $\varepsilon$. 
When the vector potential
is present, a semi-Lagrangian time splitting method has been introduced
in \cite{SL-TS,NUFFT}, with which the meshing strategy, $\Delta t=O(\varepsilon)$
and $\Delta x=O(\varepsilon)$ is needed in approximating the wave
functions. Another advantage of the time splitting methods is that if one is only interested in the physical observables, the time step size can be relaxed to $o(1)$, in other words, independently of $\varepsilon$, whereas one still needs to resolve the spatial
oscillations. 

When $\varepsilon\ll1$, several approximate methods
other than directly solving for the wave function have been proposed, such as the level
set method and the moment closure method based on the WKB analysis and the Wigner transform, see, for example, \cite{high freq waves,level set,multi-phase,reviewsemiclassical} for a general discussion.

The Gaussian beam method (or the Gaussian wave packet approach) is
another important approximate method, which allows accurate computation
around caustics and captures phase information (see \cite{Heller,popov,Gaussian propagation,EGB}).
This method reduces the full quantum dynamics to Gaussian wave packets
dynamics. Some efficient methods have been introduced to decompose general initial data into the sum of Gaussian wave packets (see \cite{EGB,FGBtrans}).
For suitable initial conditions, Gaussian beams are exact solutions to the Schr\"odinger equation with harmonic potentials. 
When the potentials are smooth (and therefore locally approximate by harmonic potentials), the Gaussian beam method is an approximate method for the solution of the Schr\"odinger equation. 
In recent years, this method has been extended for piecewise smooth
potentials (see \cite{interfaceGB}) and multi-level Schr\"odinger equations (see \cite{LuZhou1,LuZhou2}). We remark that, the Gaussian beam methods can be extended to the vector potential cases with no technical difficulties  (see \cite{HagedornV}).

However, for general potentials, the (first order) Gaussian beam method gives the approximate solution with
error $O(\varepsilon^{1/2})$, which is not practically satisfactory unless $\varepsilon$ is sufficiently small. To improve accuracy, the higher order Gaussian beam method has been introduced with error $C_{k, \varepsilon,T}\varepsilon^{k/2}$
(see \cite{HGBT,HEGB}), where $k$ stands for the approximation order and $T$ stands for the total simulation time. But still, since the modus constant in the error bound $C_{k, \varepsilon,T}$ depends on $k$ explicitly, for fixed $\varepsilon$, increasing the order of approximation in the Gaussian beam methods may  not reduce the error.
The approximation errors in high order Gaussian beam methods are especially noticeable
when $\varepsilon$ is not very small. 
It has been shown in \cite{ErrorestiGB,HagedornV}
that, when $\varepsilon=\frac{1}{100}$, in a fairly standard test example, the second
and the third order Gaussian beams give even larger error than the first
order approximation. 





Another approach to improve accuracy proposed by Faou, Gradinaru,
and Lubich in \cite{lubich Hagedorn} is to make use of the Hagedorn
wave packets. It has been shown in \cite{Hagedornraising} that
when the vector potential is absent and the scalar potential is of
a more general class, Hagedorn wave packets are asymptotic solutions
to the Schr\"odinger equation with error $O(\varepsilon^{1/2})$.
Furthermore, with higher order Hagedorn wave packets,
the error can be reduced to $O(\varepsilon^{k/2})$, where $k\in\mathbb{N}^{+}$. 
It is worth emphasizing that, higher order Hagedorn wave packets do not rely on 
cut-off functions, and can be proved to effectively reduce the approximation error for all $\varepsilon \in (0,1]$.
In \cite{lubich Hagedorn}, Faou, Gradinaru, and Lubich for
the first time turned the Hagedorn functions into a computational
tool for the semi-classical Schr\"odinger equation
with only scalar potentials. Recently, Gradinaru and Hagedorn proposed in \cite{Newsplitting} a new algorithm to solve the semi-classical Schr\"odinger
equation without vector potentials by the Hagedorn wave packets approach, which converges quadratically in the time step and linearly in $\varepsilon$. In \cite{HagedornV}, Zhou has extended this method to the vector potential cases while providing a rigorous proof for the higher order convergence with the Galerkin approximation. 

Recently, Russo and Smereka in \cite{GB Trans,GB Trans2} proposed a new method based on the so-called Gaussian wave packet transform, which reduces the quantum dynamics to Gaussian wave packet dynamics together with the time evolution of a rescaled quantity $w$, which satisfies another equation of the Schr\"odinger type, in which the modified potential becomes time dependent. 
This transform is closely related to the Gaussian wave packet method: the parameters defining the evolution of the packet satisfy a set of ordinary differential equations, so it can be accurately solved with small computational effort, while the rescaled wave function $w$ is much less oscillatory that the original wave function $\psi$, and can therefore be accurately represented with a relatively small number of degrees of freedom.
At variance with the standard Gaussian wave packet method, however,  the formulation is actually exact, which means there is no approximation error of any kind.  Furthermore, such rescaled wave function becomes even less oscillatory as the system approaches the semi-classical limit. This allows direct computation of the 
Schr\"odinger equation near the semi-classical limit even for several space dimensions.  The extension of the approach  to the multidimensional cases has been presented in \cite{GB Trans2}. 

In this paper, we aim to extend this approach to the vector potential cases. We reformulate the original Schr\"odinger equation \eqref{eq:main equation} so that it is more amenable for numerical simulation, and the Gaussian wave packet transform is carried out. In the $w$ equation obtained from the transform, a convection term shows up due to the presence of the vector potential, whose propagation velocity depends on the parameters of the Gaussian wave packet, and is thus time dependent. We then propose two first-order-in-time numerical methods to solve the $w$ equation, both of which are based on the spectral approximation and the operator splitting technique, and can naturally be extended to higher orders. The first approach is to divide the Hamiltonian into the kinetic part, the convection part and the potential part, with the first order operator splitting method. The second approach is to, instead, divide the Hamiltonian into the kinetic part and the convection-potential part. We argue that both methods can be extended to higher order ones, while in terms of efficiency, the latter one is the better choice. Rigorous stability analysis is also presented, which indicates both methods are unconditionally stable. 

 To give a heuristic presentation of this approach, the rest of the paper is organized as follows. In Section 2, a brief review for the background knowledge of the Schr\"odinger equation is presented and a brief review of the Gaussian wave packet method is given. In Section 3, we carry out the Gaussian wave packet transform in multi-dimensional cases, from which the $w$ equation is obtained. In Section 4, two first-order-in-time methods are proposed to solve the $w$ equation in one dimension, whose stability and accuracy properties are also shown. In the last section, numerous tests are performed to verify the properties of the numerical methods, where the method are made second-order in time by Strang splitting.

\section{Review of the Gaussian wave packet method}

Consider the semi-classical Schr\"odinger equation in the presence of the vector potentials,
\[
i\varepsilon\partial_{t}\psi=\frac{1}{2}\left(-i\varepsilon\nabla_{\mathbf x}-\mathbf{A}(\mathbf x)\right)^{2}\psi+V(\mathbf x)\psi.
\]
For clarity, one names  $\hat{\mathbf{P}}=-i\varepsilon\nabla_{\mathbf x}$ the canonical momentum,
and  $\hat{\mathbf{\kappa}}=\hat{\mathbf{P}}-\mathbf{A}$ is the kinetic momentum
(see \cite{FL1,FL3,Jacksonbook}).

As is discussed in the introduction, the electromagnetic field is introduced by the scalar potential $V(\mathbf x)$ and the vector potential $\mathbf{A}(\mathbf x)$. One can simplify the potential description by imposing one
more condition, namely, specifying the gauge. In fact, for any choice of a scalar function of
position $\lambda(\mathbf x)\in \mathbb{R}$, the potentials can be changed
as follows:
\begin{equation}
\mathbf{A}'=\mathbf{A}+\nabla_{\mathbf x}\lambda,\qquad V'=V.
\end{equation}

Clearly, the electric field $\mathbf{E}(\mathbf x)\in \mathbb{R}^{3}$
and the magnetic field $\mathbf{B}(\mathbf x)\in \mathbb{R}^{3}$ do not change
at all under this transformation. One natural choice is to choose
$\lambda$ so that $\nabla_{x}\cdot\mathbf{A}'=0$. This is the so-called
Coulomb gauge. In this gauge, the right hand side of the Schr\"odinger
equation (\ref{eq:main equation}) can be simplified to: 
\[
\frac{1}{2}(-i\varepsilon\nabla_{\mathbf x}-\mathbf{A})^{2}\psi+V(\mathbf x)\psi
=-\frac{\varepsilon^{2}}{2}\Delta_{\mathbf x}\psi+i\varepsilon\mathbf{A}\cdot\nabla_{\mathbf x}\psi+U(\mathbf x)\psi,
\]
where $U(\mathbf x)=\frac{1}{2}|\mathbf{A}|^{2}+V$. Actually, even if the Coulomb gauge is not chosen, one gains an extra term $\frac{i\varepsilon}{2}\nabla \cdot \mathbf A \psi$ in the equation, which is not stiff and be can also be incorporated in the $U$ term. Therefore, from this point one on, we instead consider the following form of the Schr\"odinger equation,
\begin{equation}\label{main}
i\varepsilon\partial_{t}\psi=-\frac{\varepsilon^{2}}{2}\Delta_{\mathbf x}\psi+i\varepsilon\mathbf{A}\cdot\nabla_{\mathbf x}\psi+U(\mathbf x)\psi.
\end{equation}
Note that, although the original problem is physically interesting in $3$ dimensional cases (which may reduce to $2$ dimensional problems under some circumstances), we aim to design numerical methods for this model in arbitrary dimensions. In theory, the single particle model  \eqref{eq:main equation} can naturally extend to multi-particle cases, whereas additional treatment might be needed for simulating high dimension PDE's. 
Besides, studying this model facilitates the numerical investigation of other related quantum models, such as the Pauli equation \cite{LLG-S}, which can be viewed as a vector version the Schr\"odinger equation with electromagnetic fields \eqref{eq:main equation} for spin $\frac 1 2$ particles.
 

%
Due to the intimate relations between the Gaussian wave packet based approximation methods and classical dynamics, we summarize below some basic results from classical electrodynamics. The readers may refer to, for example, \cite{Jacksonbook} for a general discussion.
Let $\bm q$ and $\bm p$ be the canonical position and momentum associated to the particle in external electromagnetic fields, and we obtain the Hamiltonian  equations of motion
\[
\dot{\bm q}=\bm p-\mathbf{A}(\bm q), \quad \dot{\bm p}=\nabla \mathbf{A}( \bm q)^{T} \bm p-\nabla U(\bm q),
\]
of the classical Hamiltonian $H(\bm q, \bm p)=\bm p^{2}-\mathbf{A}( \bm q) \cdot \bm p+U(\bm q)$.  If we denote by $\tilde {\bm q}$ and $\tilde {\bm p}$ the classical (physical) position and momentum, which are defined by the following change of variables
\[
\tilde {\bm p} =\bm p - \mathbf A (\bm q), \quad \tilde {\bm q}=\bm q, 
\]
the classical coordinates satisfy the following equations of motion (see, e.g., Chapter $12$ of \cite{Jacksonbook} for details) 
\[
\frac{d}{dt} \tilde {\bm q}= \tilde {\bm p}, \quad  \frac{d}{dt} \tilde {\bm p} = F(\tilde {\bm q}, \tilde {\bm p})= -\nabla U(\tilde {\bm p}) + \frac{d}{dt} \tilde {\bm q} \times (\nabla \times \mathbf A(\tilde {\bm p})).
\]
Here, $F$ is referred to as the Lorentz force.

For simplicity of notation, we denote the spatial coordinate with respect to the wave packet center $\bm q$ by $\bm \xi=\mathbf x-\bm q(t)$. Then 
\[
\nabla_{\mathbf x} \bm \xi={\bf 1},\quad\partial_{t} \bm \xi=-\dot{\bm q},
\]
where ${\bf 1}$ is the identity matrix.

The Gaussian beam method firstly proposed by Heller \cite{Heller} is to represent the solution by a summation of parameterized Gaussian wave packets, which take the form
\[
\varphi_{GB}=\exp\left(\frac{i}{\varepsilon} (\bm \xi \cdot \bm \alpha \bm \xi+ \bm p \cdot \bm \xi + \gamma) \right),
\] 
where $\bm \alpha$ is a complex valued symmetric matrix, $\gamma$ is a complex valued scalar, and $\bm q$ and $\bm p$ are considered as the position and the momentum of the wave packet. When the scalar potential is quadratic in $\mathbf x$ and the vector potential is linear in $\mathbf x$, it is shown in, for example, \cite{Hagedornraising,Heller} that  the wave packets are exact solutions to the  Schr\"odinger equation if the parameters satisfy the following system of equations:
\[
\left\{ \begin{array}{rcl}
 \dot {\bm q} & = &\bm p-\mathbf{A}(\bm q),
\\[2mm]
 \dot {\bm p} & = & \nabla \mathbf{A}(\bm q)^T \bm p -\nabla U(\bm q),
\\[2mm]
 \dot{ \bm \alpha }& = & -2 \bm \alpha^2 -\frac{1}{2}\nabla \nabla U(\bm q)+\nabla \mathbf{A}(\bm q)^T \bm \alpha+ \bm \alpha \nabla \mathbf{A}(\bm q)+\frac{1}{2}\nabla \nabla \mathbf{A}(\bm q) \cdot \bm p,
\\[2mm]
 \dot \gamma & = & \frac{1}{2}|\bm p|^{2}-U(\bm q)+i\varepsilon{\rm Tr}(\bm \alpha),
\end{array}\right.
\]
where ${\rm Tr (\cdot)}$ means taking the trace of a matrix.  For clarity, we present  the index form of some terms involved in the ODE systems to avoid confusion:
\[
\bigl (\nabla \mathbf{A} \bigr )_{jk}=\frac{\partial A_j}{\partial {x_k}},  \quad \bigl (\nabla \nabla U \bigr )_{jk}= \frac{\partial^2 U} {\partial {x_j}\partial {x_k}} ,\quad \bigl (\nabla \nabla \mathbf{A} \cdot \bm p )_{jk}=\sum_{\ell} \frac{\partial^2 A_{\ell} }{\partial {x_j}\partial {x_k}}\,p_{\ell}.
\]

This system of equations are well behaved as $\varepsilon \rightarrow 0$, and if properly initialized, will not develop singularities in the Riccati-type equation of $\bm \alpha$  because $\bm \alpha$ is complex valued. When the potentials are smooth, the Gaussian wave packets become approximate solutions of the  Schr\"odinger equation with error of order $O(\sqrt{\varepsilon})$ within $O(1)$ time. Besides, it is worth pointing out that, to the leading order, Heller's Gaussian wave packets are equivalent to Hagedorn's semi-classical wave packets.  

The Gaussian wave packet transform (abbreviated by GWT), which is introduced and analyzed in \cite{GB Trans, GB Trans2} for the Schr\"odinger equations with scalar potentials, is highly related to the Gaussian wave packets approach, but one of  the most significant differences is that the former one is an exact reformulation rather than an approximation. We shall provide a detailed presentation of GWT applied to the Schr\"odinger equation with vector potentials in the next section.

\section{The Gaussian wave packet transform for the  Schr\"odinger equation with vector potentials}

\subsection{The wave packet transform}

In this section, we present the Gaussian wave packet transform for the semi-classical   Schr\"odinger equation in the presence of electromagnetic field. The transform is similar to the formulations in \cite{GB Trans, GB Trans2}, but the treatment with the vector potential introduces  new challenges  and is also partially inspired by some previous work \cite{LLG-S, NUFFT,SL-TS,HagedornV}. 

We start by considering the following ansatz
\begin{equation}\label{ansatz}
\psi(\mathbf x,t)=W(\bm \xi,t)\exp \left(f(\bm \xi,t)\right):=W(\bm \xi,t)\exp \left(i\left( \bm \xi^{T} \bm \alpha_R \bm \xi +\bm p^{T}\bm \xi+\gamma_2 \right)/\e \right),
\end{equation}
where $\bm\xi=\mathbf x- \bm q$, $\bm \alpha_R$ is a real-valued symmetric matrix and $\gamma_2$ is a complex-valued scalar.
As is discussed in \cite{GB Trans}, we have to require that $\bm \alpha_R$ is real so that we can eventually arrive at a well-posed equation, which is different from other Gaussian wave packet based approaches because in those methods the Hessian matrix is assumed to be complex-valued. This {\em ansatz\/} is constructed to capture the oscillatory phase  and the Gaussian profile of the wave packet with the quadratic expansion in space with respect to the beam center, and uses the quantity $W$ to keep track of the rest of the information. And as we shall see, after making a change of variable, the equation that $W$ satisfies reduces to a non-oscillatory equation. 

We obtain, by lengthy but direct calculations (the readers may refer to Appendix A for more detailed calculations.)
\begin{equation}\label{eq:timep}
e^{-f}\partial_{t}\psi=\left(W_{t}-\nabla_{\bm \xi}W^T \dot{ \bm q}\right)+\frac{i}{\e}W\left(\bm \xi^{T}\dot{ \bm \alpha_{R}}\bm \xi+\bm \xi^{T}\dot{ \bm p}+\dot{\gamma}_{2}-2 \bm \xi^{T}\bm \alpha_{R}\dot{\bm q}-\bm p^{T}\dot{\bm q}\right),
\end{equation}
and
\begin{align}
e^{-f}\frac{i\e}{2}\Delta_{\bm x} \psi&=\frac{i\e}{2}\Delta_{\bm \xi}W- \nabla_{\bm \xi}W^{T}\left(2 \bm \alpha_{R}\bm \xi+\bm p\right)-{\rm Tr}(\bm \alpha_{R})W \nonumber \\
& \quad \quad -\frac{2i}{\e}\bm \xi^{T} \bm \alpha_R^2 \bm \xi W-\frac{2i}{\varepsilon}\bm \xi^T \bm \alpha_R \bm p W -\frac{i}{2\e}|\bm p|^2W. \label{eq:kineticp}
\end{align}

And by expanding the potentials with respect to the beam center, we have
\begin{equation}
U(\bm x)=U(\bm q)+\bm \xi^T \nabla U(\bm q)+\frac{1}{2}\bm \xi^T \nabla \nabla U(\bm q) \xi +U_{r},
\end{equation}
and
\begin{equation}
\bm A(\bm x) = \mathbf{A}(\bm q)+\nabla \mathbf{A}(\bm q)  \bm \xi+\mathbf{A}_q+\mathbf{A}_{r}
\end{equation}
where 
\[
\mathbf{A}_{(1)}=\mathbf{A}(\bm \xi+\bm q)-\mathbf{A}(\bm q)-\nabla \mathbf{A}(\bm q) \bm \xi,\quad \mathbf{A}_q= \frac{1}{2}(\bm \xi^T \nabla)^2 \mathbf{A}(\bm q),\quad \mathbf{A}_{r}=\mathbf{A}_{(1)}-\mathbf{A}_q,
\]
and
\[
U_r=U(\bm \xi+\bm q)-U(\bm q)-\bm \xi^T \nabla U(\bm q)-\frac{1}{2}\bm \xi^T \nabla \nabla U(\bm q)\bm \xi.
\]
For clarity, we represent $\mathbf{A}_q$ in the index form in the following
\[
\bigl ( \mathbf{A}_q \bigr )_{j}=\sum_{\ell,\ell'} \frac{\partial^2 A_j (\bm q)}{\partial_{x_\ell}\partial_{x_{\ell'}}}\xi_{\ell} \xi_{\ell'}.
\]

Thus, we obtain the following  (see Appendix for more details)
\begin{eqnarray*}
 e^{-f}\mathbf{A}(x)^T \nabla_{\bm x}\psi & = &  \nabla_{\bm \xi}W^T  \mathbf{A}(\bm q)+\nabla_{\bm \xi}W^T (\mathbf{A}(\bm \xi+\bm q)-\mathbf{A}(\bm q)) \\
  &   & +\frac{2i}{\e}\bm \xi^T \bm \alpha_R \mathbf{A}(\bm q)+ \frac{i}{\e}\bm \xi^T \left( \nabla \mathbf{A}(\bm q)^T\bm \alpha_R+\bm \alpha_R\nabla \mathbf{A}(\bm q)\right)\bm \xi W +\frac{2i}{\e}\mathbf{A}_{(1)}^T \bm \alpha_R \bm \xi W\\
  &   & +\frac{i}{\e} \left( \mathbf{A}(\bm q)+ \nabla \mathbf{A}(\bm q) \bm \xi + \mathbf{A}_q + \mathbf{A}_r \right)^T \bm p W.
\end{eqnarray*}
Here, we have used the fact that
\[
\bm \xi^T  \nabla \mathbf{A}(\bm q)^T\bm \alpha_R \bm \xi=\bm \xi^T \bm \alpha_R\nabla \mathbf{A}(\bm q)\bm \xi,
\]
since $\bm \alpha_R$ is symmetric.

We require the same bi-characteristic equations as in the Gaussian beam method
\[
\left\{ \begin{array}{rcl}
 \dot {\bm q} & = &\bm p-\mathbf{A}(\bm q),
\\[2mm]
 \dot {\bm p} & = & \nabla \mathbf{A}(\bm q)^T \bm p -\nabla U(\bm q),
\end{array}\right.
\]
Hence, with the Gaussian wave packet transform ansatz \eqref{ansatz}, equation \eqref{main} becomes  
\begin{eqnarray*}
W_{t} & = & -\frac{i}{\e} W \bm \xi^T \left(\dot{\bm \alpha}_{R}+2\bm \alpha_{R}^{2}+\frac{1}{2}\nabla \nabla U(\bm q)- \nabla \mathbf{A}(\bm q)^T\bm \alpha_R-\bm \alpha_R\nabla \mathbf{A}(\bm q)-\frac{1}{2}\nabla \nabla \mathbf{A}(\bm q)\cdot \bm p\right)\bm \xi\\
 &  & +W\left(\dot{\gamma}_{2}-\frac{1}{2} |\bm p|^{2}+U(\bm q)-i\varepsilon {\rm Tr} (\bm \alpha_{R})\right)\\
 &  & +\nabla_{\bm \xi}W^T \left(\mathbf{A}(\bm\xi+\bm q)-\mathbf{A}(\bm q)\right)-2\bm \xi^{T} \bm \alpha_R \nabla_{\bm \xi} W \\
 &  & + \frac{i \varepsilon}{2} \Delta_{\bm \xi} W- \frac{i}{\e} \left(U_{r}-2\bm \xi^T \bm \alpha_{R}\mathbf{A}_{(1)}-\mathbf{A}_{r}^T \bm p\right)W.
\end{eqnarray*}


Observe that, if we take
\begin{equation}\label{eq:gamma2}
\dot{\gamma}_{2}=\frac{1}{2} |\bm p|^{2}-U(\bm q)+i\varepsilon{\rm Tr}(\bm \alpha_{R}),
\end{equation}
the $W$ equation can be obviously simplified. However, we should not determine the $\bm \alpha_R$ equation directly from the $W$ equation.  As was analyzed in \cite{GB Trans, GB Trans2}, in the Gaussian wave packet transform, $\bm \alpha_R$ is considered as the real part of the complexed valued symmetric matrix $\bm \alpha$, and $\bm \alpha$ satisfies the same equation for the Hessian matrix in the standard Gaussian beam method. Namely, $\bm \alpha$ satisfies 
\begin{equation}
\dot{\bm \alpha}=-2 \bm \alpha^{2}-\frac{1}{2}\nabla \nabla U(\bm q)+\nabla \mathbf{A}(\bm q)^T\bm\alpha+\bm \alpha \nabla \mathbf{A}(\bm q)+\frac{1}{2}\nabla \nabla \mathbf{A}(\bm q) \cdot \bm p,
\end{equation}
whose real part is 
\begin{equation} \label{eq:alphaR}
\dot{\bm \alpha}_{R}=2\bm \alpha_{I}^{2}-2\bm \alpha_{R}^{2}-\frac{1}{2}\nabla\nabla U(\bm q)+\nabla \mathbf{A}(\bm q)^T \bm \alpha_R+\bm \alpha_R\nabla \mathbf{A}(\bm q)+\frac{1}{2}\nabla \nabla \mathbf{A}(\bm q) \cdot \bm p.
\end{equation} 
And, $\bm \alpha_I$, which is the imaginary part of $\bm \alpha$, is also a real-valued $N\times N$ matrix, and it satisfies
\[
\dot{\bm \alpha_I}=-2\bm \alpha_I \bm \alpha_R -2 \bm \alpha_R \bm \alpha_I + \nabla \mathbf{A}(\bm q)^T\bm  \alpha_I + \bm \alpha_I \nabla \mathbf{A}(\bm q).
\]

Then, with \eqref{eq:gamma2} and \eqref{eq:alphaR},  the $W$ equation reduces to
\begin{eqnarray*}
 W_{t} & = & \nabla_{\bm \xi}W^T \left(\mathbf{A}(\bm \xi+\bm q)-\mathbf{A}(\bm q)\right)-2 \bm \xi^{T} \bm \alpha_R \nabla_{\bm \xi} W\\
 &  &+ \frac{i \varepsilon}{2} \Delta_{\bm \xi} W- \frac{i}{\e} \left(U_{r}-2\bm \xi^T \bm \alpha_{R}\mathbf{A}_{(1)}-\mathbf{A}_{r}^T \bm p+2\bm \xi^T \bm \alpha_I^2 \bm \xi \right)W.
\end{eqnarray*}

At last, we introduce the change of variables as in \cite{GB Trans, GB Trans2}, $W(\bm \xi,t)=w(\bm \eta,t)$,
where
\begin{equation}\label{changeV}
\bm \eta= {\bm B \bm \xi}/\sqrt{\varepsilon}.
\end{equation}
Observe that, with this change of variable the $W$ function has $O(\sqrt{\e})$ support in $\xi$ while the $w$ function has $O(1)$ support in $\eta$. As we shall show in the following, due to this particular relation \eqref{changeV}, the $w$ function is not oscillatory neither in space nor in time. Here, $\bm B$ will be specified by the following equation,
\begin{equation}\label{eq:B}
\dot{\bm B}=-2\bm B \bm \alpha_R + \bm B \nabla \mathbf{A}(\bm q),
 \end{equation}
with $\bm B (0)= \sqrt{\bm \alpha_I (0)}$.  Actually, we can relate $\bm B$ to  other quantities that we have defined. By straightforward calculation, one obtains
 \[
 \frac{d}{dt}(\bm B^T \bm B)=-2\bm B^T \bm B \bm \alpha_R- 2 \bm \alpha_R \bm B^T \bm B + \nabla \mathbf{A}(\bm q)^T \bm B^T \bm B +\bm B^T\bm B \nabla \mathbf{A}(\bm q).  
 \]
 This means, since we take $\bm B(0)=\sqrt{\bm \alpha_I(0)}$, it holds that
 \begin{equation}\label{Brelation}
 \bm B(t)^T\bm B(t)= \bm \alpha_I (t), \quad \forall t \ge0.
 \end{equation}
Also, note that $\bm B(t)$ may not necessarily be symmetric even though $\bm B(0)$ is symmetric, because the equation of $\bm B(t)$ does not preserve the symmetry.  
 We remark that it has been proven in \cite{Gaussian propagation,HagedornV} that $\alpha_I$ will remain positive definite if initialized properly.  According to \eqref{Brelation}, this implies, $\bm B(t)$ will remain invertible with the initial condition given above. Therefore, $\bm \xi=\sqrt{\e} \bm B^{-1}\bm \eta$ is always well-defined.


With this change of variable, the $W$ equation becomes the following equation of $w$
\begin{align}
w_{t}=&\frac{i}{2}{\rm Tr}\left( \bm B^T \nabla_{\bm \eta} \nabla_{\bm \eta} w \bm B \right)-2i\bm \eta^{T}(\bm B^T)^{-1}\bm \alpha_I^2 \bm B^{-1}\bm \eta w  \nonumber \\
&+\frac{1}{\sqrt{\varepsilon}} \nabla_{\bm \eta}w^T \bm B \mathbf{A}_{(1)} \nonumber
\\
&+\frac{1}{i\varepsilon}\left(U_{r}-2\sqrt{\e}\mathbf{A}_{(1)}^T \bm \alpha_{R}\bm B^{-1}\bm \eta-\mathbf{A}_{r}^T \bm p\right)w.\label{weq}
\end{align}
Observe that the support of the $w$ function in the $\bm \eta$ variable is $O(1)$, 
\[
\mathbf{A}_{(1)}=O(\e),\quad \mathbf{A}_{r}=O(\e^{\frac 3 2}), \quad U_{r}=O(\e^{\frac 3 2}).
\]
Thus, we conclude that
\[
\frac{1}{\sqrt{\varepsilon}}\mathbf{A}_{(1)}=O(\sqrt{\varepsilon}),
\]
\[
\frac{1}{i\varepsilon}\left(U_{r}-2\sqrt{\e}\mathbf{A}_{(1)}^T \bm \alpha_{R}\bm B^{-1}\bm \eta-\mathbf{A}_{r}^T \bm  p\right)=O(\sqrt{\varepsilon}),
\]
so the $w$ equation is not stiff. Moreover, if one drops those $O(\sqrt{\varepsilon})$
terms, one expects to recover the leading order Gaussian beam method, as is shown in \cite{GB Trans,GB Trans2}. However, note that due to the fact the GWT parameters are time dependent, the coefficients in the $w$ equation are also time dependent, which gives rise to  additional challenges for numerical implementations. 

To summarize, with the Gaussian wave packet transform, the  Schr\"odinger equation \eqref{main} equivalently transforms to the $w$ equation \eqref{weq} together with the ODE system of Gaussian wave packet parameters
\begin{equation} \label{psys}
\left\{ \begin{array}{rcl}
 \dot {\bm q} & = & \bm p-\mathbf{A}(\bm q),\\[2mm]
  \dot {\bm p} & = & \nabla \mathbf{A}(\bm q)^T \bm p -\nabla U(\bm q),\\[2mm]
   \dot {\bm \alpha} & = & -2 \bm \alpha^2 -\frac{1}{2}\nabla \nabla  U(\bm q)+\nabla \mathbf{A}(\bm q)^T\bm \alpha+\bm \alpha\nabla \mathbf{A}(\bm q)+\frac{1}{2}\nabla \nabla \mathbf{A}(\bm q) \cdot \bm p,
\\[2mm]
 \dot \gamma_2 & = & \frac{1}{2}|\bm p|^{2}-U(\bm q)+i\varepsilon{\rm Tr}(\bm \alpha_{R}),
\\[2mm]
  \dot{\bm B} & = & -2\bm B \bm \alpha_R + \bm B \nabla \mathbf{A}(\bm q).
\end{array}\right.
\end{equation}

Compared with the Gaussian beam method, the equations of the parameters $\bm q$, $\bm p$ and $\bm \alpha$ are exactly the same, but the equations for $\gamma$ and $\gamma_2$ differ slightly,
\[
\frac{d}{dt}(\gamma-\gamma_2)=-\e {\rm Tr}(\bm \alpha_I).
\]
However, it's worth emphasizing that the Gaussian wave packet transform is an exact reformulation, rather than an approximation like the Gaussian beam method.

We remark however that as $\e$ becomes smaller and smaller, then better and better accuracy is required in the solution of the ODE system for the Gaussian wave packet parameters, since the phase of the wave function $\psi$ contains a factor $1/\e$, and small error in the phase may be amplified by a large factor.

\subsection{Initial condition and extensions}

In this part, we summarize the procedures of the Gaussian wave packet transform method  and discuss the possible initial conditions that this approach can handle. We start with an initial condition  in a Gaussian profile:
\begin{equation}\label{init1}
\psi(\mathbf x,0)=\exp \left( i \left( \bm \xi_0^T \bm C \bm \xi_0+\bm \xi_0^T \bm p_0 + \delta \right)/\e \right),
\end{equation}
where $\bm \xi_0=\mathbf x -\bm x_0$, $\bm C$ is a symmetric complex valued matrix, with its imaginary part $\bm C_I$ positive definite and $\delta$ is a complex valued scalar.

 Clearly, it follows from ansatz \eqref{ansatz} and the change of variable \eqref{changeV} that, if we specify $\bm B(0)=\sqrt{\bm C_I}$, the initial condition for the $w$ equation is 
 \begin{equation}\label{initw}
 w(\bm \eta,0)=e^{-\bm\eta^T \bm \eta}.
 \end{equation}
Note that $\bm C_I$ is positive definite, so that $\bm B(0)$ is well defined. Besides, we want to emphasize that the initial condition for $w$ is actually independent of the Gaussian wave packet parameters. Now, we are ready to give the initial conditions for the parameters:
\[
\bm q(0)=\bm x_0,\quad \bm p(0)=\bm p_0, \quad \bm \alpha (0)=\bm  C, \quad \gamma_2(0)=\delta, \quad \bm B(0)=\sqrt{\bm C_I}.
\]


Therefore, given the Schr\"odinger equation \eqref{main} with initial condition \eqref{init1}, one can alternatively solve the $w$ equation \eqref{weq} and the ODE system \eqref{psys} with initial conditions defined above, and the solution to the original problem is reconstructed by
\begin{equation}
\psi (\mathbf x,t)=w(\bm B\bm \xi/\sqrt \e,t)\exp \left( i \left(\bm \xi^{T} \bm \alpha_R \bm \xi +\bm p^{T}\bm \xi+\gamma_2 \right)/ \e \right), \label{recon}
 \end{equation}
 where we recall that $\bm \xi=\mathbf x - \bm q(t)$. We remark that, not only the $w$ equation \eqref{weq} and the system \eqref{psys} contain no stiffness as $\e \rightarrow 0$, but also the high oscillation is removed from the initial condition for the $w$ equation. Therefore, meshes which do not depend on $\e $ are allowed for the $w$ equation for accurate numerical approximations.  
 
 However, in practice, the initial conditions may not take the form of a Gaussian wave packet. In case of more general initial functions, one solution is to decompose the function into a sum of Gaussian wave packets, and solve the problem by the Gaussian wave packet transform for each wave packet individually. Due to the linearity of the problem, the sum of all the reconstructed solutions is exactly the solution to the original problem. The decomposition itself is an active area of research. To our knowledge, one of the best results is given by Qian and Ying in \cite{FGBtrans}, where they proposed a fast algorithm to decompose a wide class of smooth functions into a superposition of Gaussian wave packets.  
 
 On the other hand, due to the flexibility of this transform, as discussed in \cite{GB Trans,GB Trans2}, the Gaussian wave packet transform if capable of handing more general initial conditions in natural way. Now, consider the initial conditions of the following form
 \begin{equation}\label{init2}
 	\psi(\mathbf x,0)=f(\mathbf x-\bm q_0) \exp \left( i \left( g(\mathbf x-\bm q_0)+ (\mathbf x -\bm q)^T \bm p_0 + \delta \right)/\e \right),
 \end{equation} 
where $f$ and $g$ are complex valued smooth functions, with $f$ bounded. Without loss of generality, we demand $g(0)=0$, $\nabla g (0)= \mathbf 0$, $\im\, g$ is convex, and if we define $\bm C=\bm C_R+i \bm C_I= \nabla \nabla g(0)/2$, we also require that $\bm C_I$ is positive definite. Now the initial condition for the $w$ equations becomes
\[
	w(\bm \eta,0)=f(\bm \xi) \exp \left( i\left(g(\bm \xi)- \bm \xi^T\bm C_R\bm \xi \right)/\e \right),
\]
where $\bm \xi=\mathbf x- \bm q_0=\sqrt \e \bm B^{-1}\eta$. And the initial conditions for the system \eqref{psys} are 
\[
	\bm q(0)=\bm x_0,\quad \bm p(0)=\bm p_0, \quad \bm \alpha (0)= \bm C, \quad \gamma_2(0)=\delta, \quad \bm B(0)=\sqrt{\bm C_I}.
\]

From this point of view, the Gaussian wave packet transform can deal with initial conditions of a wider class, which implies one does not have to necessarily decompose a smooth initial condition into a superposition of Gaussian wave packets. We would like to explore in the future whether a better decomposition method can be introduced so that the whole method will become more efficient. 

We remark that, as is pointed out in \cite{GB Trans,GB Trans2}, the Gaussian wave packet transform facilitates efficient calculation of a large family of physical observables without reconstructing the wave function. In other words, those physical observables can be expressed in terms of the $w$ function and the Gaussian wave packet parameters. In the presence of the vector potentials, the expressions can be derived without additional difficulties. For examples, the expectation value of the position can be obtained by
\begin{equation} \label{xave}
<\bm x>= \int \bm x |\psi|^2 d \bm x = \bm q(t)+ \frac{\e^{\frac{N+1}{2}}}{|\bm B|} e^{-{2\im (\gamma_2)}/{\e}}B^{-1} \int \bm \eta |w|^2 d \bm \eta.
\end{equation} 
Here, $|\bm B|$ denotes the determinant of $\bm B$.

We would like to conclude this section with comments on the computational domain of the $w$ equation. If the initial condition to the Schr\"odinger equation \eqref{main} is of a Gaussian profile, its effective width is $O(\sqrt{\e})$. And due to the relations between the $\mathbf x$ variable and the $\eta$ variable, we learn that we should truncate the $\eta$ space to a $O(1)$ domain in order to enforce the periodic boundary condition  with negligible domain truncation error. When the initial condition is of a general type, clearly for a wide class of initial conditions, it suffices to prescribe a $O(1)$ domain for the $w$ equation.

\section{Numerical implementation}

In the section, we propose and analyze numerical methods to implement the Gaussian wave packet transform approach. The numerical simulation of the ODE system \eqref{psys} is standard. Besides, since the system does not depend on the $w$ equation \eqref{weq}, one can numerically solve the parameters till any given time with arbitrary time steps and to arbitrary accuracy. Note that, the time steps in solving the ODE system has to be $\e$ dependent, since the numerical error will be magnified by a factor of $\e^{-1}$. However, this is not considered as a challenge due to the existence of various higher ODE solvers. For example, if a forth order ODE solver is applied to the system, the time steps in numerically integrating the ODE system can be taken as $o(\e^{-1/4})$.

Because in general, solving the ODE system is much cheaper than numerically solving the partial differential equations, we would assume the the parameters are solved firstly with minimal error and we, therefore, only need to focus on numerical methods to the $w$ equation.

\subsection{Description of the numerical method for the $w$ equation}

In this part, we give a full construction of the numerical method for the $w$ equation. We start with time splitting of the $w$ equation.
Here, we propose two ways to split the Hamiltonian.

\paragraph{Three step scheme}
The first way agrees with what is done in \cite{SL-TS} for the  Schr\"odinger equation: for every time step $t\in[t_{n},t_{n+1}]$,
one solves the {\bf kinetic step}
\begin{equation}
w_{t}=\frac{i}{2}{\rm Tr}\left( \bm B^T \nabla_{\bm \eta}  \nabla_{\bm \eta}  w \bm B \right),\quad t\in[t_{n},t_{n+1}];\label{eq:split1}
\end{equation}
followed by the {\bf potential step}
\begin{equation}
w_{t}=-2i \bm \eta^{T}(\bm B^T)^{-1}\bm \alpha_I^2\bm  B^{-1}\bm \eta w+\frac{1}{i\varepsilon}\left(U_{r}-2\sqrt{\e}\mathbf{A}_{(1)}^T  \bm \alpha_{R}\bm B^{-1}\eta-\mathbf{A}_{r}^T\bm p\right)w,\quad t\in[t_{n},t_{n+1}],\label{eq:split2}
\end{equation}
and followed by the {\bf convection step}
\begin{equation}
w_{t}=\frac{1}{\sqrt{\varepsilon}} \nabla_{\bm \eta}w^T \bm  B \mathbf{A}_{(1)},\quad t\in[t_{n},t_{n+1}].\label{eq:split3}
\end{equation}

\paragraph{Two step scheme}
The second way is to combine the potential step and the convection step together. Namely, for every time step $t\in[t_{n},t_{n+1}]$, one solves the {\bf kinetic step}
\begin{equation}
w_{t}=\frac{i}{2}{\rm Tr}\left( \bm B^T \nabla_{\bm \eta}  \nabla_{\bm \eta}  w \bm B \right),\quad t\in[t_{n},t_{n+1}];\label{eq:split21}
\end{equation}
followed by the {\bf convection-potential step}
\begin{equation} \label{eq:split22}
w_{t}=-2i\eta^{T}(\bm B^T)^{-1}\bm \alpha_I^2\bm  B^{-1}\bm \eta w+\frac{1}{i\varepsilon}\left(U_{r}-2\sqrt{\e}\mathbf{A}_{(1)}^T \bm \alpha_{R}\bm B^{-1}\bm \eta-\mathbf{A}_{r}^T\bm p\right)w
\end{equation}
\begin{equation*}
+\frac{1}{\sqrt{\varepsilon}} \nabla_{\bm \eta}w^T\bm  B \mathbf{A}_{(1)},\quad t\in[t_{n},t_{n+1}].
\end{equation*}


Compared to the  Schr\"odinger equation, the $w$ equation is significantly cheaper to solve because the $w$ function is essentially not oscillatory. However, even with operator splitting, there is no practical way to solve the $w$ equation by analytical solutions in each substep because the coefficients are time dependent.  

For simplicity, we present the numerical method to solve the $w$ equation \eqref{weq} in one
dimension with periodic boundary condition. The extension to multidimensional
cases is straightforward.
We assume, on the computation domain $[a,b]$, a uniform spatial grid
$\eta_{j}=a+j\Delta \eta$, $j=0,\cdots N-1$, where $N=2^{n_{0}}$, $n_{0}$
is an positive integer and $\Delta \eta=\frac{b-a}{N}$. We also assume
uniform time steps $t_{n}=n\Delta t$, $n=0,\cdots,n_{\rm max}$. The numerical approximation of the $w$ function at $t=t_n$ is denoted by {$w^n$}, with components denoted by {$w^n_j$}. 
The construction
of numerical methods is based on the operator splitting technique. For clarity, we only discuss the first order splitting for this moment. The extension to higher order splitting methods will be discussed later. 

\paragraph{Kinetic step}
{The kinetic step can be very effectively solved in Fourier space.}
In the one dimensional case, we define the Fourier coefficients of $w^k$ in the following way
\[
\zeta_l=\frac{2\pi l}{b-a},\quad \hat w^k_l=\sum_{j=0}^{N-1}w_j^k e^{-i {\zeta}_l (\eta_j-a)}, \quad l=-\frac{N}{2},\cdots, \frac{N}{2}-1.
\]
By applying the Fourier transform to the  { kinetic step}, we obtain 
\[
\hat w_t =   - \frac{i}{2}  B^2 \zeta^2 \hat w.
\]
Thus, the analytical solution to the kinetic step is
\begin{equation}\label{kstep}
\hat w (\zeta,t)=\exp \left( -i  \int_{t_0}^t \frac 1 2 B^2\zeta^2 d s \right) \hat w (\zeta,t_0).
\end{equation}
Recall that, we have assumed $B(t)$ has been solved in advance with great accuracy. Then to numerically solve $w(t)$ with a specific order of accuracy in time, one just needs to apply some quadrature rules of the corresponding order to approximate the time integral in \eqref{kstep}. Note that, in the time approximation with quadratures, the spatial variable $\zeta$ acts as a parameter.  We denote the approximation of the time integral at $\zeta=\zeta_l$ by {$K^n_l$}, where
\[
{	
	K^n_l =\int_{t_n}^{t_{n+1}}  \frac{1}{2}
	B^2 \zeta_l^2 ds
	}
\]
and obviously $|\exp(-i K_l)|=1$. Then, the numerical solution to the kinetic step is given by the following:
\[
	w^{k+1}_j=\frac{1}{N}\sum_{l=-N/2}^{N/2-1} e^{ {-}iK_l}\hat w^k_l e^{i \zeta_l (\eta_j-a)},\quad j=0,\cdots,N-1.
\]

\paragraph{Potential step}
Next, to simplify the notations, we introduce for the potential step and the convection step
 \[
	F(\eta,t)=-2i\eta  B^{-1}  \alpha_I^2 B^{-1}\eta+\frac{1}{i\varepsilon}\left(U_{r}-2\sqrt{\e} {A}_{(1)}  \alpha_{R} \bm B^{-1}\eta-{A}_{r} p\right),
\]
and 
\[
	G(\eta,t)=\frac{1}{\sqrt{\varepsilon}}  B {A}_{(1)},
\]
where $G(\eta,t)$ is real-valued and $F(\eta,t)$ is purely imaginary. Thus, the potential step becomes
\begin{equation} \label{eq:pots}
w_t = F(\eta,t) w,
\end{equation}
and the convection step becomes
\begin{equation}\label{eq:cons}
w_t =G(\eta,t) w_{\eta}.
\end{equation}

For the potential step \eqref{eq:pots}, the analytical solution is
\[
	w(\eta,t)=\exp \left( \int_{t_0}^t F(\eta,s)\,ds \right)w(\eta,t_0).
\]
 And to obtain a numerical scheme of certain order, one just need to apply a corresponding quadrature rule to approximate the time integral. We denote the approximation of the integral at $\eta=\eta_j$ by $F^k_j$, i.e.
\[
	F^n_j \approx \int_{t_n}^{t_{n+1}} F(\eta_j,s)\,ds
\]
and obviously $|\exp(F_j)|=1$. Then, the numerical method to the potential step is given by
\[
	{w^{n+1}_j=e^{F^n_j}w^n_j, \quad j=0,\cdots, N-1.}
\] 

\paragraph{Convection step}
The convection part \eqref{eq:cons}, however, has to be treated differently because on the discrete spatial grids, there is no analytical solution that one can make use of. Here, we use the semi-Lagrangian method as proposed in \cite{SL-TS} to solve the convection part. This method consists of two parts: backward characteristic tracing and interpolation. To
compute the data $w_{j}^{n+1}$ which are approximations of $w(t_{n+1},\eta_j)$, we firstly trace backwards along
the characteristic line, along which $w$ remains constant:
\begin{equation}
\frac{d\eta(t)}{dt}=-G(\eta(t),t),\quad \eta(t_{n+1})=\eta_{j},\label{eq:char}
\end{equation}
for time interval $[t_{n},t_{n+1}]$.  If we denote $\eta(t_{n})=\eta_{j}^{0}$,
obtained by numerically solving the ODE (\ref{eq:char}) backwards
in time,  then by the method of characteristics, $w(\eta_j,t_{n+1})=w(\eta_j^0,t_n)$. However, the numerical approximations of  $w(\eta_j^0,t_n)$ are not in general known since $\eta_j^0$ may not be the grid points. Therefore, certain interpolation is needed to calculate $w_{j}^{n+1}\approx w(\eta_j,t_{n+1})=w(\eta_j^0,t_n)$. The semi-Lagrangian method for the convection step with the polynomial interpolation has been studied in \cite{SL-TS}, and the one with the spectral interpolation is studied in \cite{NUFFT}. We remark that, with the Nonuniform FFT algorithms, the numerical cost of the spectral interpolation in the convection step can be reduced to $O(N \log N)$ (see \cite{NUFFT,nufft6}), which is comparable to the numerical cost in the kenetic step. In this paper, we choose to apply the spectral interpolation to compute $\{w_{j}^{n+1}\}$ based on $\{\eta_j^0\}$, which is spectrally accurate in space and unconditionally stable.

In summary, for the convection step, we proposed a semi-Lagrangian method which consists of numerically solving the characteristic equations backward in time and applying spectral interpolation to compute the point values.

With the full description of the scheme above, we name the whole scheme the three-part semi-Lagrangian time splitting spectral method (abbreviated by SL-TS3). 

\paragraph{Convection-Potential step}
Next, we discuss the numerical approximation for the convection-potential step:
\begin{equation}
w_t=G(\eta,t)w_{\eta}+ F(\eta,t)w.
\end{equation}
Following the characteristics,
\[
\frac{d\eta(t)}{dt}=-G(\eta(t),t),
\]
we rewrite the convection-potential step as,
\[
w_t(\eta(t),t)=F(\eta(t),t)w(\eta(t),t).
\]
So the exact solution to this step is
\[
w(t,\eta(t))=\exp\left(\int_{t_0}^t F(\eta(s),s)ds\right)w(t_0,\eta(t_0)).
\]

In order to make use of this solution, one can numerically approximate the time integral 
\[
\int_{t_0}^t F(\eta(s),s)ds \approx  \Pi_m F,
\]
where $\Pi_m F$ is an order $m$ approximation with corresponding quadrature points. For example, one may take 
\[
\Pi_1 F = (t-t_0) F(\eta(t_0),t_0)\quad {\rm and}\quad \Pi_2 F = \frac{1}{2}(t-t_0)\left( F(\eta(t_0),t_0)+ F(\eta(t),t) \right).
\] 
Also, note that $\Pi_m F$ is different from $F_j$ defined before because here $\eta(t)$ is no longer constant in time, but it is a time dependent trajectory. 

Obviously, if $\eta(t)$ is exactly known, $\Pi_m F$ is indeed an order $m$ approximation of that time integral. However, note that $\eta(t)$ is yet to be approximated as well with a certain order of accuracy, so direct discretization $\Pi_m F$ with approximate $\eta(t)$ may not yield an order $m$ accuracy. But, as we show in the following, the error in the backward tracing of $\eta (t)$ is much smaller than the numerical approximation error of the time integral of $F$.  

To simplify the analysis, we define $H(t)=\int_{t_n}^t F(\eta(s),s)ds $, and consider the system     
\[
\begin{cases}
\displaystyle \frac{d}{dt}H(t)=F(\eta(t),t), & H(t_n)=0,  \\[3mm]
\displaystyle \frac{d}{dt}\eta(t)=-G(\eta(t),t), & \eta(t_{n+1})=\eta_j. \\
\end{cases}
\]
And we are solving for $\eta(t_n)$ and $H(t_{n+1})$, which we call the backward-forward step. 

The key observation is, clearly, $F = O(1)$ but $G=O(\sqrt{\e})$. Also, the equation of $\eta(t)$ is independent of the equation for $H(t)$.  Therefore, we can first apply an ODE solver for the $\eta(t)$ equation with the numerical error reduced by  a factor  $\sqrt{\e}$, and subsequently solve for $H(t)$ with the same or a different ODE solver. And the accuracy order of $H(t)$ follows by standard numerical analysis. Note that, even if $G=O(1)$, the accuracy argument is still valid, but the fact that $G=O(\sqrt{\e})$ makes the error in finding $\eta(t)$ much smaller than that in $H(t)$. Besides, it is worth emphasizing that, the backward-forward step has only $O(N)$ cost, which is much less than the cost of the kinetic step, which is $O(N \log N)$. 


We remark that, we would introduce the convection-potential step for the $w$ equation mainly due to the following two reasons:
\begin{enumerate}
\item  The convection velocity $-G=O(\sqrt \eps)$, and the $w$ function is essentially not oscillatory, which is very different from the original Schr\"odinger equation. Therefore, the exact solution to the convection-potential step can easily be approximated numerically with satisfactory accuracy.
\item  If one aims for a higher order method in time for the $w$ equation, a higher order operator splitting technique needs to be applied. In this sense, a two-part splitting is more efficient than a three-part splitting.  We shall elaborate this point in Section \ref{sec:comment}.
\end{enumerate}

With the full description of the whole scheme above, we name this whole scheme the two-part semi-Lagrangian time splitting spectral method (abbreviated by SL-TS2).

\subsection{Stability and Accuracy}

In this section, we will show that the SL-TS2 and the SL-TS3 methods for the one dimensional $w$ equation are unconditionally stable when the error in solving the backward tracing step is negligible. The extension to the multidimensional cases is straightforward, and in Section \ref{sec:num} we numerically test the performance of these methods in some multidimensional problems.

\subsubsection{Stability}

For the stability analysis, because the $w$ equation is linear,  in the perspective of the operator splitting, it suffices to show the stability for each sub-step.

We start by studying the convection step. As is analyzed in \cite{NUFFT, SL-TS,suli}, if the error in backward characteristic tracing is negligible, the semi-Lagrangian method for the convection part is unconditionally stable. In the $w$ equation,  the convection velocity $G(\eta,t)$ is time-dependent due its dependence on the GWT parameters, so one should solve the characteristic equations for every time step. Thanks to the fact that the convection velocity $-G=O(\sqrt \eps)$, it is affordable to approximate the backward characteristics with minimal error. 

Define $w^n$ to be the numerical approximation of the $w$ function at the beginning of  the convection step at $t=t_n$ with $w^n_j$ as components. {And we denote by $w^n_I$ the spectral approximation of the $w$ function based on $w^n$, namely
\[
w^{n}_I(x)=\frac{1}{N}\sum_{l=-N/2}^{N/2-1} \hat w^n_l e^{i \zeta_l (x-a)}, \quad a<x<b,
\]
where we recall $\hat w^k_l$ are the Fourier coefficients based on $w^n_j$.
One immediate result is, with the periodic boundary condition, 
\[
	{\left\Vert w^n \right\Vert_{l^2}=\left\Vert w_I^n \right\Vert_{L^2}=\left\Vert w_I^n (\eta_j) \right\Vert_{l^2}.}
\]
Here, the norms used are defined in the following
\[
\left\Vert w^n \right\Vert_{l^2}= \left( \frac{b-a}{N} \sum_{j=1}^N |w^n_j|^2 \right)^{\frac 1 2}, \quad \left\Vert w_I^n \right\Vert_{L^2}= \left( \int_a^b |w^n_I (x)|^2 d x \right)^{\frac 1 2}.
\]

{According to the stability analysis in \cite{NUFFT, SL-TS,suli},} 
if $\eta_j^0$ are computed with negligible error, 
\[
	\left\Vert \tilde W^{n+1} \right\Vert_{l^2}= \left\Vert w^{n}_I (\eta_j^0) \right\Vert_{l^2} \leqslant  (1+C\Delta t)   \left\Vert w^{n}_I (\eta_j) \right\Vert_{l^2} =(1+C\Delta t) \left\Vert w^n \right\Vert_{l^2},
\]
where the constant $C$ can be taken as $\max_t \|\partial_{\eta}G(\eta,t)\|_{L^{\infty}}$. In particular, $C$ is independent of $\Delta t$, and $C=0$ when $G(\eta,t)$ is divergence free.  Here, the inequality follows from the fact that the norm of the projection operator on the pseudo-spectral subspace is no larger than $1$ and  by a standard estimate of the analytical evolution operator of the convection equation, see \cite{NUFFT,suli}.  

Next, in the potential step, recalling that $F_j$ are purely imaginary, we conclude that,
\[
|w^{k+1}_j|=|e^{F_j}w^k_j|=|w^k_j|,
\]
which implies $\Vert w^{k+1}\Vert_{l^2}=\Vert w^k \Vert_{l^2}$. Hence, the numerical method for the potential step is unconditionally stable.

Based on the results in the convection step and the potential step, we can easily show that the numerical method for the convection-potential step in SL-TS2 is also unconditionally stable.

Finally, we aim to prove that the numerical method for the kinetic step is unconditionally stable. The calculation is similar to that in \cite{TSSP}, but to our best knowledge, the stability analysis with time dependent coefficient has not been considered before. The proof is given by the following calculations.
\begin{align*}
\frac{1}{b-a}\Vert w^{k+1} \Vert_{l^2} & =  \frac{1}{N}\sum_{j=0}^{N-1} |w^{k+1}_j|^2 \\
 & =  \frac{1}{N}\sum_{j=0}^{N-1} \left|\frac{1}{N}\sum_{l=-N/2}^{N/2-1} e^{-iK_l}\hat w^k_l e^{i \zeta_l (\eta_j-a)}\right|^2 \\
 & =  \frac{1}{N^3} \sum_{j=0}^{N-1} \sum_{l_1=-N/2}^{N/2-1} \sum_{l_2=-N/2}^{N/2-1}e^{-i(K_{l_1}-K_{l_2})} \hat w^k_{l_1} \overline{\hat W^k_{l_2}}e^{i(\eta_j-a)(\zeta_{l_1}-\zeta_{l_2})} \\
 & =  \frac{1}{N^3} \sum_{j=0}^{N-1} \sum_{l_1=-N/2}^{N/2-1} \sum_{l_2=-N/2}^{N/2-1}e^{-i(K_{l_1}-K_{l_2})} \hat w^k_{l_1} \overline{\hat w^k_{l_2}}e^{i {2 \pi j(l_1-l_2)}/{N}} \\
 & =  \frac{1}{N^2}  \sum_{l=-N/2}^{N/2-1} |\hat w^k_{l}|^2 ,
 \end{align*}
 where $\hat w^k_{l}$ denotes the Fourier coeffecients. Then, we obtain
\begin{align*}
\frac{1}{b-a}\Vert w^{k+1} \Vert_{l^2}  & =  \frac{1}{N^2}  \sum_{l=-N/2}^{N/2-1} \left|\sum_{j=0}^{N-1}w_j^k e^{-i \zeta_l (\eta_j-a)}\right|^2 \\
 & =  \frac{1}{N} \sum_{j=0}^{N-1}|w_j^k|^2 =\frac{1}{b-a}\Vert w^{k} \Vert^{l^2}.
\end{align*}
Here, we have used the identities
\[
\sum_{j=0}^{N-1} \exp\left(i\frac{2\pi(k-l)j}{N}\right)=\sum_{l=-N/2}^{N/2-1} \exp\left(i \frac{2\pi(k-l)j}{N}\right)= N\delta_{kl},
\]
where $\delta_{kl}=1$ when $k-l=mN$, $m\in\Z$, otherwise $\delta_{kl}=0$. We remark that, this result suggests that no matter how accurate $K_l$ are in approximating the phase multiplier in the kinetic step, the spectral method  preserves the $l^2$ norm exactly.

To sum up, we have shown that both the SL-TS2 method and the SL-TS3 method for the $w$ equation are unconditionally stable, which facilitates us to study the convergence results.

\subsubsection{Accuracy and other comments} \label{sec:comment}

The numerical error of a scheme based on the the Gaussian wave packet transform is mainly due to the following three sources: numerical approximation for the Gaussian wave packet parameters, the operator splitting for the $w$ equation, the numerical approximation for the $w$ function in each sub step.
The errors introduced by the numerical approximation of the Gaussian wave packet parameters and by operator splitting are standard. 

For the operator splitting error, if one aims for a high order method, compared with SL-TS3, the SL-TS2 method would be more efficient because as is analyzed in \cite{H3P}, it is harder and takes more sub-steps to construct a high order operator splitting methods with three parts. For example, to our best knowledge, the fourth order operators splitting method with two parts requires at least 5 steps, while the fuorth order operator splitting method with three parts requires at least 12 steps (see \cite{Horder,H3P}).      

Next, we turn to discuss the error in solving the $w$ equation for each sub-step. Observe that the smoothnesses of the vector potential ${\bf A}$ and the scalar potential $V(x)$ indicate that $F(\eta,t)$, $G(\eta,t)$ and their derivatives are of order $O(1)$ and order $O(\sqrt{\e})$, respectively. In addition, recall the $w$ function is not oscillatory. By standard arguments, one concludes that the method has spectral accuracy in space for the potential and kinetic step. If one chooses to use spectral interpolation for the semi-Lagrangian method in the convection part, the method also has  spectral accuracy in space in the convection or convection-potential step. Also, one safely concludes that the accuracy in time for each step is dominated by the quadrature rules used and the ODE solver chosen for the characteristic equations.  

 To conclude this section, we want to emphasize that, the domain  of the $w$ equation is chosen in an \emph{ad} \emph{hoc} way, which may introduce some noticeable error as well. This issue has been discussed in \cite{GB Trans, GB Trans2}, and in numerical tests we numerically check whether the computational domain has been chosen large enough. 

\subsection{Reference solution and error evaluation}

As discussed in the introduction, one of the alternative approach to solve the Schr${\rm \ddot o}$dinger with vector potentials is the semi-Lagrangian time-splitting method introduced in \cite{NUFFT, SL-TS} (abreviated by SL-TS). This method is based on the Fourier spectral method for the potential and the kinetic step, and a semi-Lagrangian method in the convection step. The method is unconditionally stable, has spectral accuracy in space if a spectral interpolation method is used in the convection step, and can be constructed to be arbitrary order accurate in time. However, in the semi-classical regime, due to the $O(\e)$ small scale oscillation in space and time, one has to take the following meshing strategy to get the accurate approximation of the wave function:
\[
\Delta t = O(\e),\quad \quad \Delta x=O(\e).
\]
In the numerical tests of time convergence, we will use this method with sufficiently fine mesh to compute the reference solution.

Note that, by the Gaussian wave packet transform approach, the computation domain is fixed in the $\eta$ space. However, the corresponding domain in the $x$ space is varying in time due to the relation:
\[
\mathbf x=\bm q(t)+\sqrt{\e} \bm B(t)^{-1}\bm \eta.
 \] 
Therefore, we need to apply certain interpolation to match the data points with the reference solution. Here, we choose to use spectral interpolation to maintain high accuracy.

Besides, note that the corresponding domain in the $x$ space may not be square. Then in one dimensional tests, it is still reasonable to compute the absolute error, but in high dimensional cases, it is more convenient to compute the relative error instead.  The readers may refer to \cite{GB Trans2} to see more discussions on this issue.   

\section{Numerical examples} \label{sec:num}

In this section,  we aim to test the SL-TS2 method and the SL-TS3 method proposed in the previous sections. We would like to remind the readers that, the numbers in the abbreviations indicate how many parts the Hamiltonian is divided into. We choose to use the Strang splitting in time and spectral interpolation in the convection step. Hence, we expect to verify the second order convergence in time and the spectral convergence in space. Also, we want to test whether the proposed methods can handle more general initial conditions. Both one dimensional and two dimensional tests are provided.
%

\subsection*{Example 1} In this problem, we consider the one dimensional Schr${\rm \ddot o}$dinger equation with the scalar potential and the vector potential given by 
\[
V(x)=1+\cos(x),\quad A(x)=\sin(x).
\]
The initial condition is given by a Gaussian wave packet
\begin{equation}\label{initc}
\psi_0(x)=\left( \frac{2 {\rm Im}(\alpha_0)}{\pi \eps} \right)^{1/4} \exp \left(i \frac{\alpha_0 (x-q_0)^2}{\e}+i\frac{p_0(x-q_0)}{\e} +i\frac{\gamma_0}{\e} \right) ,
\end{equation} 
where
\[
q_0=\frac{\pi}{4},\quad p_0=-\frac{1}{2},\quad \alpha_0=i,\quad \gamma_0=0.
\]
Note that, the initial wave function has been properly normalized.

Before investigating the convergence behavior in time and in space separately, we want to compare the respective errors from time discretization and from spatial discretization. To this end, we fix $\e= \frac{1}{256}$, and compare the solutions of the GWT based SL-TS3 method with various $\Delta \eta$ and $\Delta t$ to the reference solution at $t=0.5$. The reference solution is computed by the SL-TS method with highly resolved mesh: $\Delta x=\frac{2\pi \eps}{32}$, $\Delta t=\frac{\eps}{32}$. And $L^2$ errors are plotted in Figure \ref{dt_deta}, from which clearly see that when the spatial mesh is highly coarse, the numerical error is dominated by the spatial error, which are manifested by the "plateaus" the error curves reach when  $\Delta t$ decreases for $\Delta \eta= \frac{\pi}{4}$, $\Delta \eta=\frac{\pi}{5}$ and $\frac{\pi}{6}$. But, when the spatial mesh is fine enough (e.g. $\Delta \eta \le \frac {\pi}{8}$), the numerical error is dominated by the time discretization.  And we can clearly observe second order convergence in time for $\Delta \eta= \frac{\pi}{8}$, $\Delta \eta=\frac{\pi}{9}$ and $\frac{\pi}{10}$. Similar tests have been done with the GWT based SL-TS2 method, but the results are rather similar, so we choose to omit them in this work. We carry out detailed convergence studies in the following.

\begin{figure} 
\begin{centering}
\includegraphics[scale=0.75]{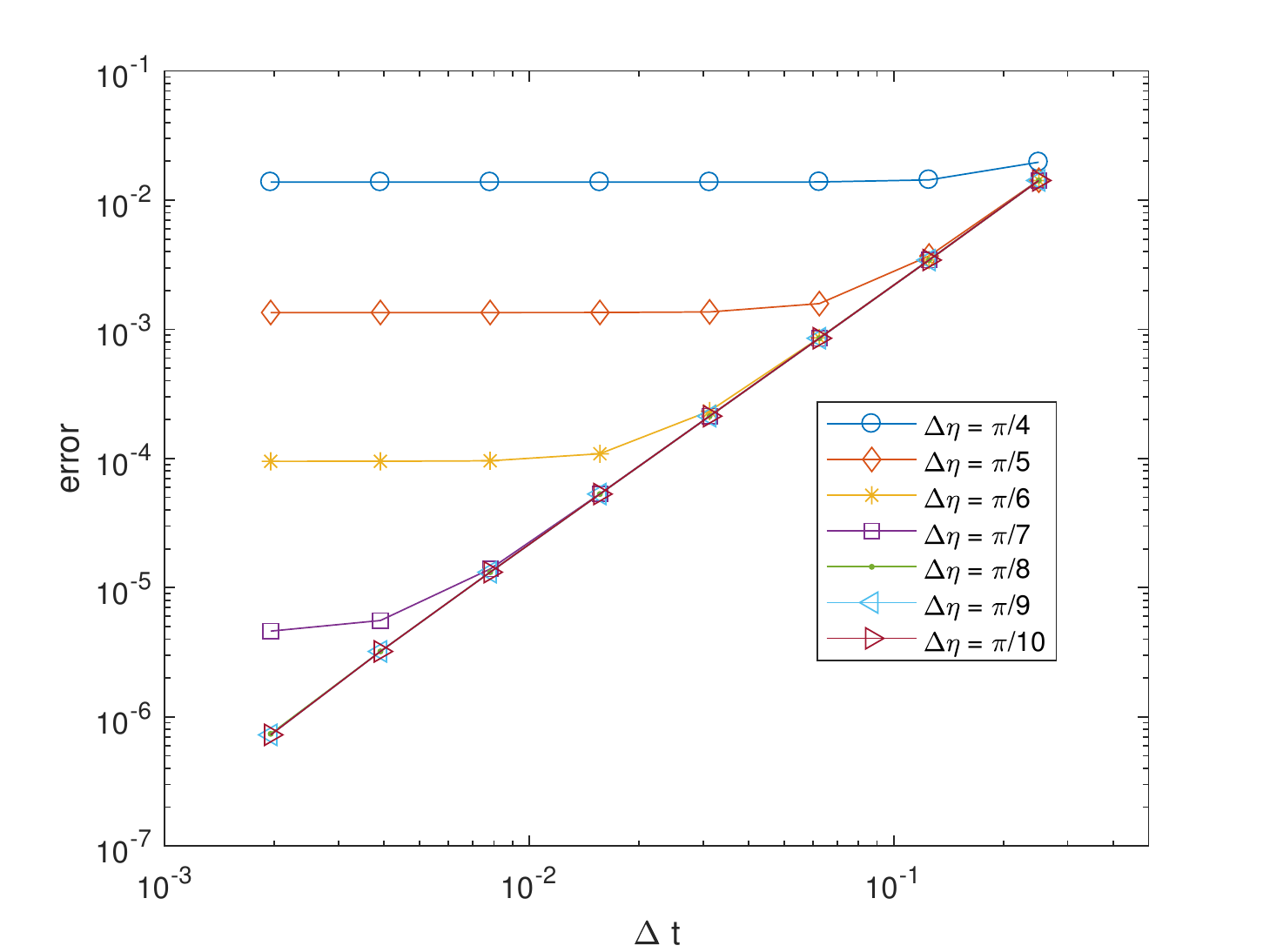} \par
\end{centering}
\caption{(Example 1) $\varepsilon=\frac{1}{256}$. Reference solution SL-TS: $\Delta x=\frac{2\pi \eps}{32}$, $\Delta t=\frac{\eps}{32}$.  GWT based SL-TS3 method: $\Delta \eta=\frac{\pi}{4}$, $\Delta \eta=\frac{\pi}{5}$, $\Delta \eta=\frac{\pi}{6}$, $\Delta \eta=\frac{\pi}{7}$, $\Delta \eta=\frac{\pi}{8}$, $\Delta \eta=\frac{\pi}{9}$ and $\Delta \eta=\frac{\pi}{10}$;  $\Delta t=\frac{1}{4}$, $\frac{1}{8}$, $\frac{1}{16}$, $\frac{1}{32}$, $\frac{1}{64}$, $\frac{1}{128}$, $\frac{1}{256}$ and $\frac{1}{512}$.}
\label{dt_deta}
\end{figure}

In this part, we would like to test the convergence in time. The reference solution is computed on $[-\pi,\pi],$ till the comparison time $T=0.5$. For various $\eps$, the corresponding meshing strategy is
$\Delta x=\frac{2\pi \eps}{32}$ and $\Delta t=\frac{\eps}{32}$,
which is sufficiently refined to generate benchmark solutions.  

To test the convergence in time of the SL-TS3 method, we apply the fourth order Runge-Kutta method to the ODE system of the parameters with time step $\delta t$, and for the $w$ equation,  we choose the computation domain to be $[-2\pi,2\pi]$ with well-resolved spatial mesh $\Delta \eta=4\pi/1024$. For various $\eps=\frac{1}{256}$, $\frac{1}{512}$, $\frac{1}{1024}$, $\frac{1} {2048}$, and various time steps $\Delta t=\frac{1}{8}$, $\frac{1}{16}$, $\frac{1}{32}$, $\frac{1}{64}$, $\frac{1}{128}$, we calculate the error of the numerical solution in $L^2$ norm. Note that, since with resolved spatial mesh, the error in spatial discretization is minimal compared to the error in time. Also, since the time step $\delta t$ for the ODE system is independent of the choice of $\Delta t$, and  the ODE's are cheaper to solve, we choose $\delta t=\Delta t/40$, which makes the error in evolving the parameters negligible. Recall that, the ODE's have to be computed very accurately since the error in the phase function is magnified by $\varepsilon^{-1}$. Also, we remark that other accurate ODE integrators  can also be applied here, such as symplectic integrators. Thus, one expect the numerical error is dominated by the error in time discretization. We plot the results in Table \ref{table:test1}, which clearly shows the second order convergence in $\Delta t$. Besides, we also notices that the error seems to be independent of $\eps$, which is expected as well because the Gaussian wave packet transform is an exact reformulation. Next, we repeat the tests with the SL-TS2 method, and plot the results in Table \ref{table:test2}. We see clearly, the errors with the SL-TS3 method is nearly the same as the previous case. 


\begin{table}
\small
  \centering
  \begin{tabular}{ c|c | c| c|c|c} \hline
    $L^2$ error &  $\Delta t= \frac{1}{8}$ & $ \Delta t= \frac{1}{16}$ & $\Delta t= \frac{1}{32}$ & $\Delta t= \frac{1}{64}$ & $\Delta t=  \frac{1}{128}$   \\ \hline
$\varepsilon= \frac{1}{256}$ & 3.431e-3 & 8.457e-4 & 2.056e-4 & 4.653e-5  & 1.047e-5
    \\ \hline
 $\varepsilon=\frac{1}{512}$ & 3.428e-3 & 8.467e-4 & 2.071e-4  & 4.701e-5 & 1.056e-5
    \\ \hline
$\varepsilon= \frac{1}{1024}$ & 3.444e-3 & 8.491e-4 & 2.064e-4 & 4.673e-5  & 1.035e-5
    \\ \hline
 $\varepsilon=\frac{1}{2048}$ & 3.447e-3 & 8.496e-4 & 2.066e-4  & 4.676e-5 & 1.034e-5
    \\ \hline
\end{tabular}
\caption{(Example 1) $\varepsilon=\frac{1}{256}$, $\frac{1}{512}$, $\frac{1}{1024}$, $\frac{1}{2048}$. Reference solution SL-TS: $\Delta x=\frac{2\pi \eps}{32}$, $\Delta t=\frac{\eps}{32}$. Comparing with GWT based SL-TS3 method: $\Delta \eta=\frac{4\pi}{1024}$,  $\Delta t=\frac{1}{8}$, $\frac{1}{16}$, $\frac{1}{32}$, $\frac{1}{64}$, $\frac{1}{128}$.}
  \label{table:test1}
\end{table}


\begin{table}
\small
  \centering
  \begin{tabular}{ c|c | c| c|c|c} \hline
    $L^2$ error &  $\Delta t= \frac{1}{8}$ & $ \Delta t= \frac{1}{16}$ & $\Delta t= \frac{1}{32}$ & $\Delta t= \frac{1}{64}$ & $\Delta t=  \frac{1}{128}$   \\ \hline
$\varepsilon= \frac{1}{256}$ & 3.425e-3 & 8.424e-4 & 2.049e-4 & 4.740e-5  & 1.022e-5
    \\ \hline
 $\varepsilon=\frac{1}{512}$ & 3.437e-3 & 8.464e-4 & 2.058e-4  & 4.708e-5 & 1.128e-5
    \\ \hline
$\varepsilon= \frac{1}{1024}$ & 3.443e-3 & 8.483e-4 & 2.063e-4 & 4.694e-5  & 1.081e-5
    \\ \hline
 $\varepsilon=\frac{1}{2048}$ & 3.446e-3 & 8.492e-4 & 2.065e-4  & 4.689e-5 & 1.056e-5
    \\ \hline
\end{tabular}
\caption{(Example 1) $\varepsilon=\frac{1}{256}$, $\frac{1}{512}$, $\frac{1}{1024}$, $\frac{1}{2048}$. Reference solution SL-TS: $\Delta x=\frac{2\pi \eps}{32}$, $\Delta t=\frac{\eps}{32}$. Comparing with GWT based SL-TS2 method: $\Delta \eta=\frac{4\pi}{1024}$,  $\Delta t=\frac{1}{8}$, $\frac{1}{16}$, $\frac{1}{32}$, $\frac{1}{64}$, $\frac{1}{128}$.}
  \label{table:test2}
\end{table}


We now move on to test the spatial convergence. To eliminate the error in time discretization in comparison, we use the numerical solution by the SL-TS3 method with sufficiently fine mesh as the reference solution:
\[
\Delta \eta=\frac{2\pi}{4096}, \quad \Delta  t =\frac{1}{2048}.
\] 

For $\e=\frac{1}{1024}$, $\frac{1}{2048}$, we compute the numerical solution by the SL-TS3 method with same time step as the reference solution, but with various spatial mesh size $\Delta \eta= \frac{\pi}{1}$, $\frac{\pi}{2}$, $\frac{\pi} {4}$, $\frac{\pi} {8}$ and $\frac{\pi} {16}$. The results are plotted in Table \ref{table:test3}, from which we see clearly that as $\Delta \eta$ decreases, the numerical error decays exponentially fast until it becomes minimal. We remark that, the same tests have been done with the SL-TS2 method, and the results are almost the same, so we omit this part in the paper.

\begin{table}
\small
  \centering
  \begin{tabular}{ c|c | c| c|c|c} \hline
    $L^2$ error &  $\Delta \eta= \frac{\pi}{1}$ & $ \Delta \eta= \frac{\pi}{2}$ & $\Delta \eta= \frac{\pi}{4}$ & $\Delta \eta= \frac{\pi}{8}$ & $\Delta \eta=  \frac{\pi}{16}$   \\ \hline
$\varepsilon= \frac{1}{1024}$ & 8.446e-1 & 2.955e-1 & 1.269e-2 & 7.234e-8  & 7.956e-11
    \\ \hline
 $\varepsilon=\frac{1}{2048}$ & 8.427e-1 & 2.888e-1 & 1.226e-2  & 6.371e-8 & 7.953e-11
    \\ \hline
\end{tabular}
\caption{(Example 1) $\varepsilon = \frac{1}{1024}$, $\frac{1}{2048}$. Reference solution, GWT based  SL-TS2 with $\Delta \eta=\frac{2\pi}{4096}$ and $\Delta  t =\frac{1}{2048}$. Comparing with GWT based SL-TS2 method: $\Delta t=\frac{1}{1024}$,  $\Delta \eta=\frac{\pi}{1}$, $\frac{\pi}{2}$, $\frac{\pi}{4}$, $\frac{\pi}{8}$ and  $\frac{\pi}{16}$.}
  \label{table:test3}
\end{table}

To further demonstrate the efficiency of the Gaussian wave packets transform in terms of spatial meshes, we focus on the case when $\e=\frac{1}{1024}$, and numerical solution is obtained by the SL-TS3 method with $\Delta t=\frac{1}{2048}$ and $\Delta \eta= \frac{\pi}{16}$, but the wave function $\psi(x)$ is reconstructed by the Fourier interpolation in the $w$ function on a refined mesh with $\Delta x = \frac{\pi}{8192}$, which is compared with the reference solution computed by the SL-TS method for the Schr${\rm \ddot o}$dinger with highly resolved time steps and spatial meshes, see Figure \ref{fig:comp}. Besides, we plot the computed grid values of $W(\xi-q)$ (we have shifted $W$ by the beam center to facilitate comparison) in the same plot and $w(\eta)$ in the next plot to show that, the computational grids do not resolve the wave function $\psi(x)$ at all, it only resolves the non-oscillatory $w(\eta)$, which is sufficient to ensure accurate reconstruction of the wave function $\psi(x)$. 
\begin{figure}
\begin{centering}
\includegraphics[scale=0.45]{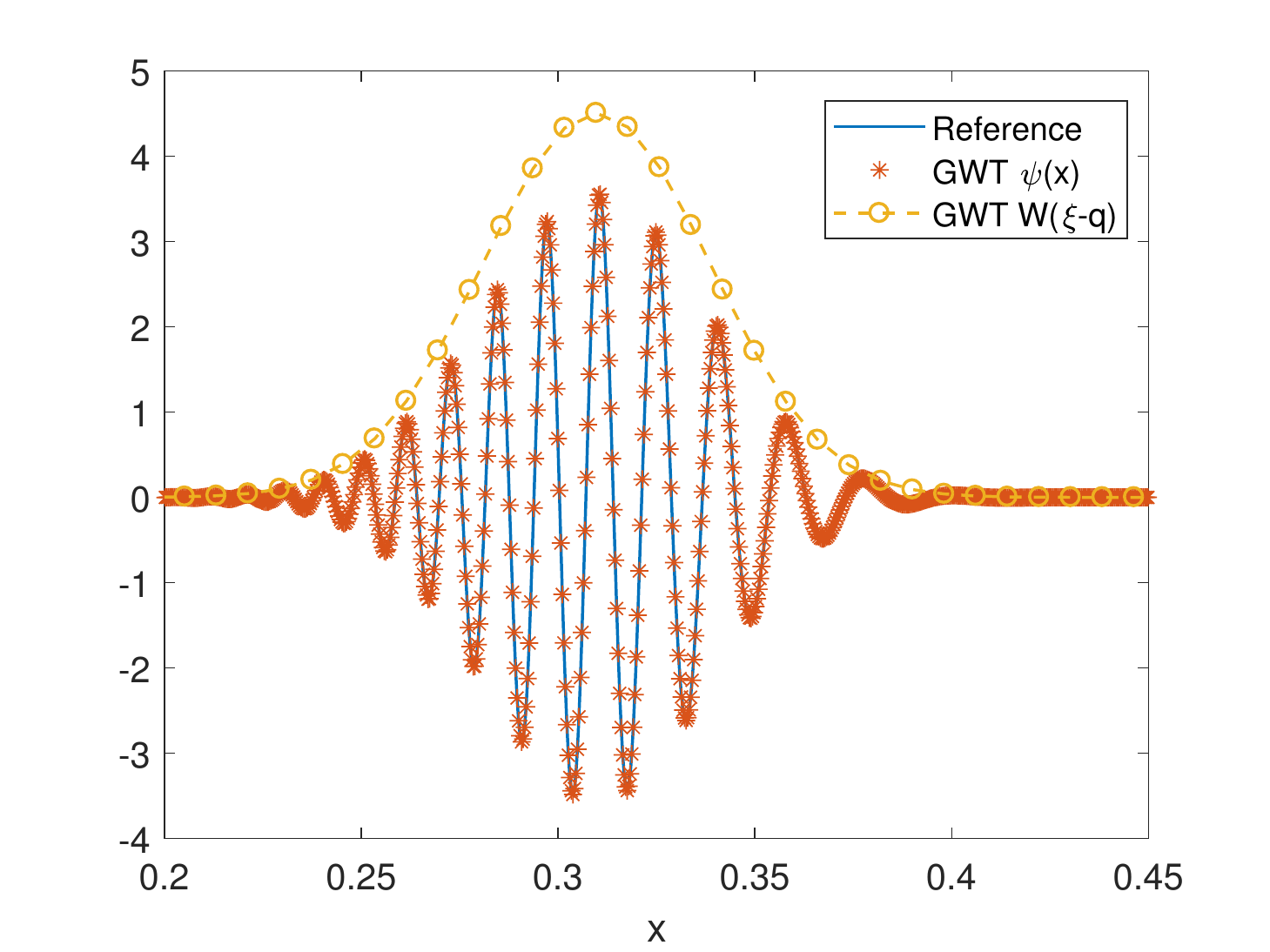} \includegraphics[scale=0.45]{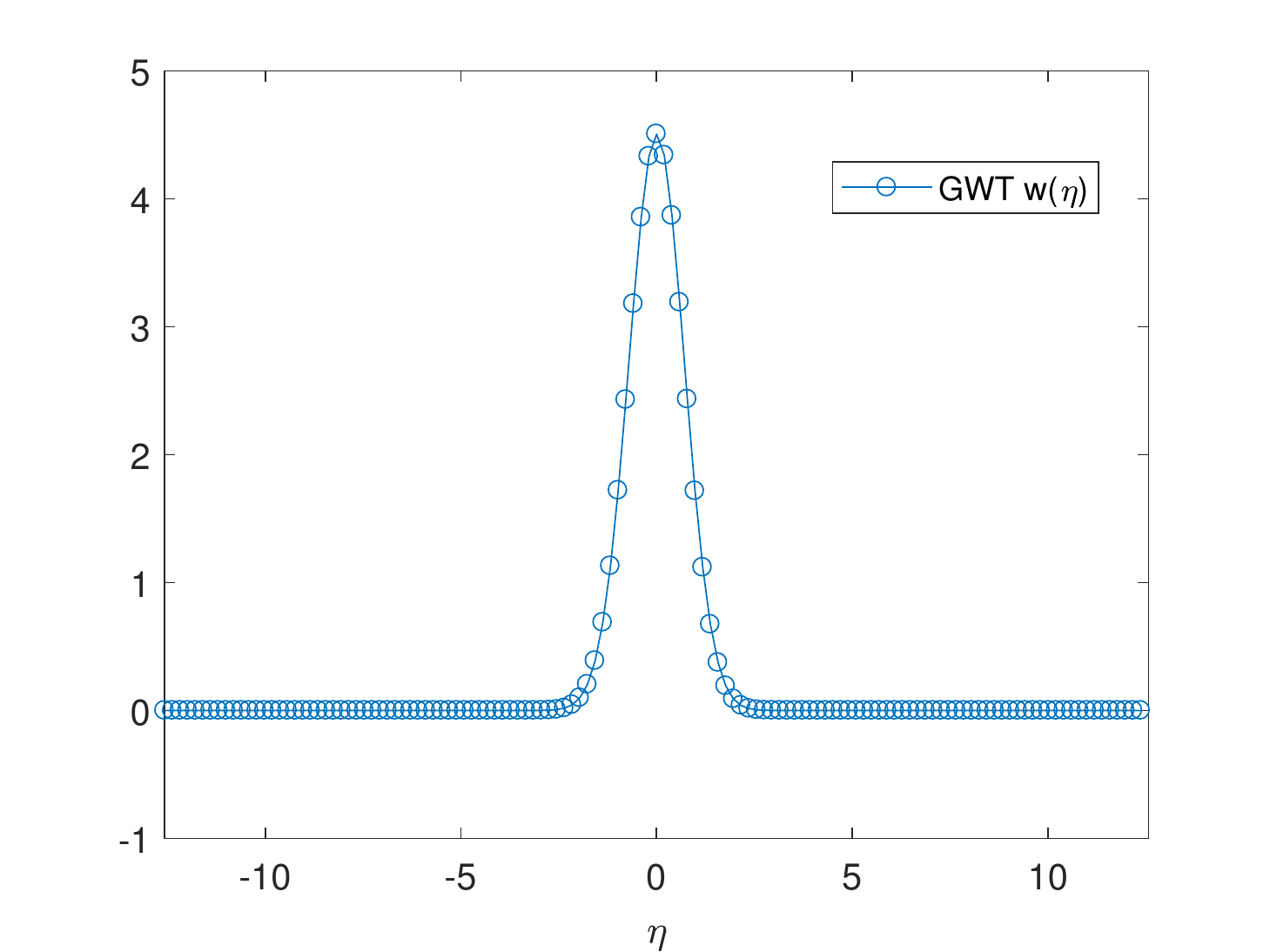} \par
\end{centering}
\caption{(Example 1) $\e=\frac{1}{1024}$.  Left: the blue line denotes the real part of reference solution $\psi(x)$, the red stars denote the real part of the reconstructed wave function $\psi(x)$ by the GWT based SL-TS3 method, the yellow circles denote the real part of the corresponding computed grid values of $W(\xi-q)$. Right: the blue circles denote the real part of the computed grid values of $w(\eta)$ with $\Delta \eta= \frac{\pi}{16}$.}
\label{fig:comp}
\end{figure}

%
  


Finally, we report the error growth of the GWT method in time in the presence of vector potentials. In general, the numerical error depends on several factors including, $\Delta t$, $\Delta \eta$, the computation domain and the semi-classical parameter $\eps$. In the work, we only test the GWT based SL-TS3 method for various $\epsilon$. To single out the effect of the vector potential, we change the potentials to
\[
V(x)=2x^2,\quad A(x)=e^{-2x^2},
\]
which means the scalar potential is harmonic and the vector potential effectively vanishes at the boundary of the computational domain.  The initial condition is still given by by a Gaussian wave packet as in \eqref{initc}, with 
\[
q_0=0,\quad p_0=0,\quad \alpha_0=i,\quad \gamma_0=0.
\]
We choose $\Delta t = \frac{1}{100}$, $\Delta \eta = \frac{\pi}{64}$ and the computational domain for $w$ equation is $[-16, \, 16]$, and we plot the errors versus time for $\varepsilon = \frac{1}{64}$, $\frac{1}{256}$ and $\frac{1}{1024}$ in Figure \ref{fig:time}. We observe that for each $\eps$, the error grows exponentially in time, and the accuracy of the GWT method is better for smaller $\eps$. This results is basically in line with the numerical studies in \cite{GB Trans, GB Trans2}, although the present of the vector potential make the performance of the GWT method slightly worse in long time simulations.

\begin{figure} 
\begin{centering}
\includegraphics[scale=0.75]{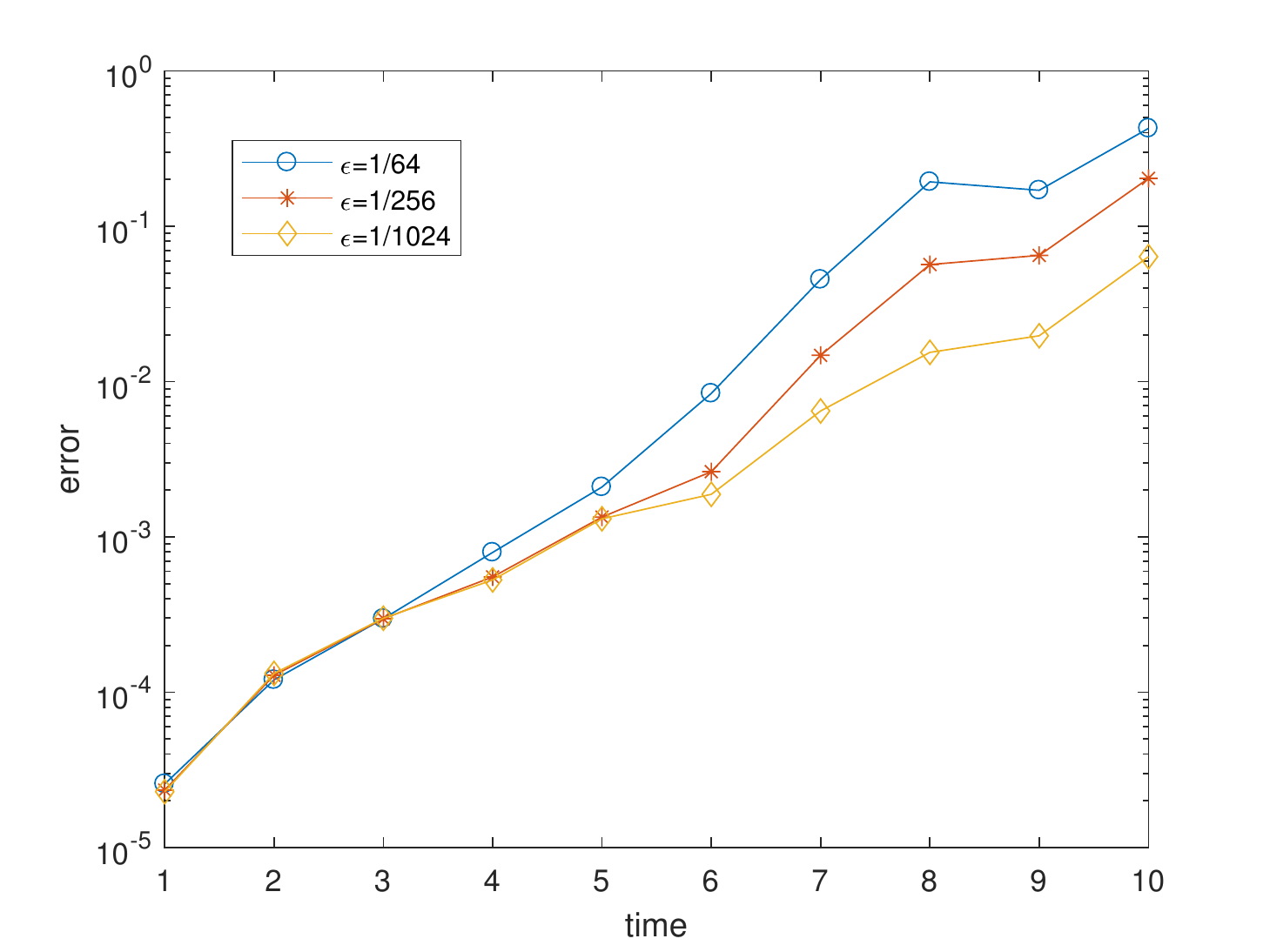} \par
\end{centering}
\caption{(Example 1) Numerical errors versus time for $\varepsilon = \frac{1}{64}$, $\frac{1}{256}$ and $\frac{1}{1024}$. $V(x)=2x^2$ and $A(x)=e^{-2x^2}$. The GWT based SL-TS3 method: $\Delta t = \frac{1}{100}$, $\Delta \eta = \frac{\pi}{64}$, computation domain $[-16, \, 16]$. Reference solutions: the SL-TS method with $\Delta t= \frac{\eps}{16}$ and $\Delta x= \frac{2\pi \eps}{32}$. }
\label{fig:time}
\end{figure}

\subsection*{Example 2} In this problem, we consider the one dimensional Schr${\rm \ddot o}$dinger with the scalar potential and the vector potential given by 
\[
V(x)=1+\cos(x),\quad A(x)=\sin(x),
\]
which is the same as the previously one. But the initial wave function is chosen as 
\[
\psi_0(x)=a \exp\left(-\frac{(x-q_0)^2+(x-q_0)^4}{\e} \right) \exp\left(i\frac{\cos(x-q_0)}{\e}\right),
\]
where $q_0=\frac{\pi}{4}$, and $a$ is chosen so that $\psi_0(x)$ is normalized. Clearly, this initial wave function is not of a Gaussian profile, but it satisfies the form \eqref{init2}, where
\[
f(y)=a,\quad g(y)=i (y^2+y^4)+ \cos(y)-1,\quad p_0=0,\quad \gamma=1.
\]
In this part, we aim to repeat the time convergence and spatial convergence test with the more general initial condition.  Again, the reference solution is computed on $[-\pi,\pi],$ till the comparison time $T=0.5$. For various $\eps$, the corresponding meshing strategy is
\[
\Delta x=\frac{2\pi \eps}{32},\quad \Delta t=\frac{\eps}{32},
\]
which is sufficiently refined to generate benchmark solutions.  

To test the convergence in time of the SL-TS3 method and the SL-TS2 method, we apply the fourth order Runge-Kutta method to the ODE system of the parameters with time step $\delta t$, and for the $w$ equation,  we choose the computation domain to be $[-2\pi,2\pi]$ with well-resolved spatial mesh $\Delta \eta=4\pi/1024$. For various $\eps$ and various time steps, we calculate the error of the numerical solution in $L^2$ norm.  We plot the results of the SL-TS3 method in Table \ref{table:test4} and those of the SL-TS2 method in Table \ref{table:test5}, which clearly show the second order convergence in $\Delta t$, although the initial condition is no longer of a Gaussian profile.

\begin{table}
\small
  \centering
  \begin{tabular}{ c|c | c| c|c|c} \hline
    $L^2$ error &  $\Delta t= \frac{1}{8}$ & $ \Delta t= \frac{1}{16}$ & $\Delta t= \frac{1}{32}$ & $\Delta t= \frac{1}{64}$ & $\Delta t=  \frac{1}{128}$   \\ \hline
$\varepsilon= \frac{1}{256}$ & 4.627e-3 & 1.155e-3 & 2.918e-4 & 7.633e-5  & 2.266e-5
    \\ \hline
 $\varepsilon=\frac{1}{512}$ & 4.619e-3 & 1.153e-3 & 2.914e-4  & 7.630e-5 & 2.267e-5
    \\ \hline
$\varepsilon= \frac{1}{1024}$ & 4.615e-3 & 1.152e-3 & 2.912e-4 & 7.629e-5  & 2.268e-5
    \\ \hline
 $\varepsilon=\frac{1}{2048}$ & 4.613e-3 & 1.152e-3 & 2.911e-4  & 7.629e-5 & 2.269e-5
    \\ \hline
\end{tabular}
\caption{(Example 2) $\varepsilon=\frac{1}{256}$, $\frac{1}{512}$, $\frac{1}{1024}$, $\frac{1}{2048}$. Reference solution SL-TS: $\Delta x=\frac{2\pi \eps}{32}$, $\Delta t=\frac{\eps}{32}$. Comparing with GWT based SL-TS3 method: $\Delta \eta=\frac{4\pi}{1024}$,  $\Delta t=\frac{1}{8}$, $\frac{1}{16}$, $\frac{1}{32}$, $\frac{1}{64}$, $\frac{1}{128}$.}
  \label{table:test4}
\end{table}


\begin{table}
\small
  \centering
  \begin{tabular}{ c|c | c| c|c|c} \hline
    $L^2$ error &  $\Delta t= \frac{1}{8}$ & $ \Delta t= \frac{1}{16}$ & $\Delta t= \frac{1}{32}$ & $\Delta t= \frac{1}{64}$ & $\Delta t=  \frac{1}{128}$   \\ \hline
$\varepsilon= \frac{1}{256}$ & 4.627e-3 & 1.158e-3 & 2.950e-4 & 7.859e-5  & 2.434e-5
    \\ \hline
 $\varepsilon=\frac{1}{512}$ & 4.619e-3 & 1.155e-3 & 2.930e-4  & 7.762e-5 & 2.353e-5
    \\ \hline
$\varepsilon= \frac{1}{1024}$ & 4.615e-3 & 1.153e-4 & 2.921e-4 & 7.695e-5  & 2.311e-5
    \\ \hline
 $\varepsilon=\frac{1}{2048}$ & 4.613e-3 & 1.152e-4 & 2.915e-4  & 7.662e-5 & 2.290e-5
    \\ \hline
\end{tabular}
\caption{(Example 2) $\varepsilon=\frac{1}{256}$, $\frac{1}{512}$, $\frac{1}{1024}$, $\frac{1}{2048}$. Reference solution SL-TS: $\Delta x=\frac{2\pi \eps}{32}$, $\Delta t=\frac{\eps}{32}$. Comparing with GWT based SL-TS2 method: $\Delta \eta=\frac{4\pi}{1024}$,  $\Delta t=\frac{1}{8}$, $\frac{1}{16}$, $\frac{1}{32}$, $\frac{1}{64}$, $\frac{1}{128}$.}
  \label{table:test5}
\end{table}


Then, we test the spatial convergence. Again, to eliminate the error in time discretization in comparison, we use the numerical solution by the SL-TS3 method with sufficiently fine mesh as the reference solution. For $\e=\frac{1}{1024}$, $\frac{1}{2048}$, we compute the numerical solution by the SL-TS3 method with the same time step as the reference solution, but with various spatial mesh size. The results are summarized in Table \ref{table:test6}, from which we see clearly that as $\Delta \eta$ decreases, the numerical error decays exponentially fast until it becomes minimal. We remark that, the same tests have been done with the SL-TS2 method, and the results are almost the same, so we omit this part in the paper.

\begin{table}
\small
  \centering
  \begin{tabular}{ c|c | c| c|c|c} \hline
    $L^2$ error &  $\Delta \eta = \frac{\pi}{1}$ & $ \Delta \eta= \frac{\pi}{2}$ & $\Delta \eta= \frac{\pi}{4}$ & $\Delta \eta= \frac{\pi}{8}$ & $\Delta \eta =  \frac{\pi}{16}$   \\ \hline
$\varepsilon= \frac{1}{1024}$ & 9.392e-1 & 3.434e-1 & 1.282e-2 & 5.744e-8  & 8.079e-11
    \\ \hline
 $\varepsilon=\frac{1}{2048}$ & 7.874e-1 & 2.888e-1 & 1.075e-2  & 4.646e-8 & 6.793e-11
    \\ \hline
\end{tabular}
\caption{(Example 2) $\varepsilon = \frac{1}{1024}$, $\frac{1}{2048}$. Reference solution, GWT based  SL-TS2 with $\Delta \eta=\frac{2\pi}{4096}$ and $\Delta  t =\frac{1}{2048}$. Comparing with GWT based SL-TS2 method: $\Delta t=\frac{1}{1024}$,  $\Delta \eta=\frac{\pi}{1}$, $\frac{\pi}{2}$, $\frac{\pi}{4}$, $\frac{\pi}{8}$ and  $\frac{\pi}{16}$.}
  \label{table:test6}
\end{table}




\subsection*{Example 3}  In this test,  we consider the three dimensional problem demonstrated in  Appendix \ref{3dex}, where the whole problem naturally decomposes into a two dimensional problem with vector potentials in $(x,y)$ and a quantum harmonic oscillator in $z$. Although the problem is three dimensional, it suffices to apply the GWT method to solve for the the two dimensional wave function $u(x,y,t)$.

The initial condition is chosen as a Gaussian wave packet in $(x,y)$ multiplied by $\phi_0(z)$,
\[
\psi_0({\bf x})=\tilde \psi_0(x,y)\phi_0(z)=a\exp\left( i\frac{\xi_0 \cdot C \xi_0+ p_0\cdot \xi_0 + \gamma_0}{\eps} \right)\phi_0(z),
\]
where $\xi=(x,y)^T-q_0$, $a$ is a normalization constant.  And for the potentials, we take
\[
A_1(y)=\sin(y), \quad A_2(x)= \sin(x), \quad \tilde V(x,y)= \cos(x)+\cos(y),
\]
and the corresponding magnetic field is given by
\[
\bm{B} = (0,0,\cos(x)-\cos(y))^T.
\] 

The reference solution is computed on $[-\pi,\pi]^2$, till the comparison time $T=1$. For various $\eps$, the corresponding meshing strategy is
\[
\Delta x=\Delta y=\frac{\pi \eps}{16},\quad \Delta t=\frac{\eps}{16},
\]
which is sufficiently refined to generate benchmark solutions.

To test the convergence in time of the SL-TS3 method, we choose Gaussian wave packet parameters as
\[
q_0=\left(0.5,0\right)^T,\quad p_0=(-2,0)^T,\quad C=i \,{\mathrm I}_2,\quad \gamma_0=0,
\]
where ${\mathrm I}_2$ is the $2\times 2$ identity matrix.
We apply the fourth order Runge-Kutta method to the ODE system of the parameters with time step $\delta t$, and for the $w$ equation,  we choose the computation domain to be $[-8,8]^2$ with well-resolved spatial mesh $\Delta \eta=8/512$. For various $\eps=\frac{1}{16}$, $\frac{1}{32}$ and $\frac{1}{64}$, and various time steps $\Delta t=\frac{1}{8}$, $\frac{1}{16}$, $\frac{1}{32}$, $\frac{1}{64}$, we calculate the relative error of the numerical solution in $L^2$ norm at time $t=0.5$ (more details in calculating the relative error can be found in \cite{GB Trans2}). We first plot the (real parts of) numerical solutions in Figure \ref{fig:2d}. To make the comparison convenient, we have shifted the $W$ by the beam center $\bm q$, and we observe that the $W$ function is much smoother than the oscillatory wave function $\psi$, but its support is $O(\sqrt{\e})$, and clearly, the $w(\bm \eta)$ is smooth with $O(1)$ sized support. 
\begin{figure}
\begin{centering}
\includegraphics[scale=0.45]{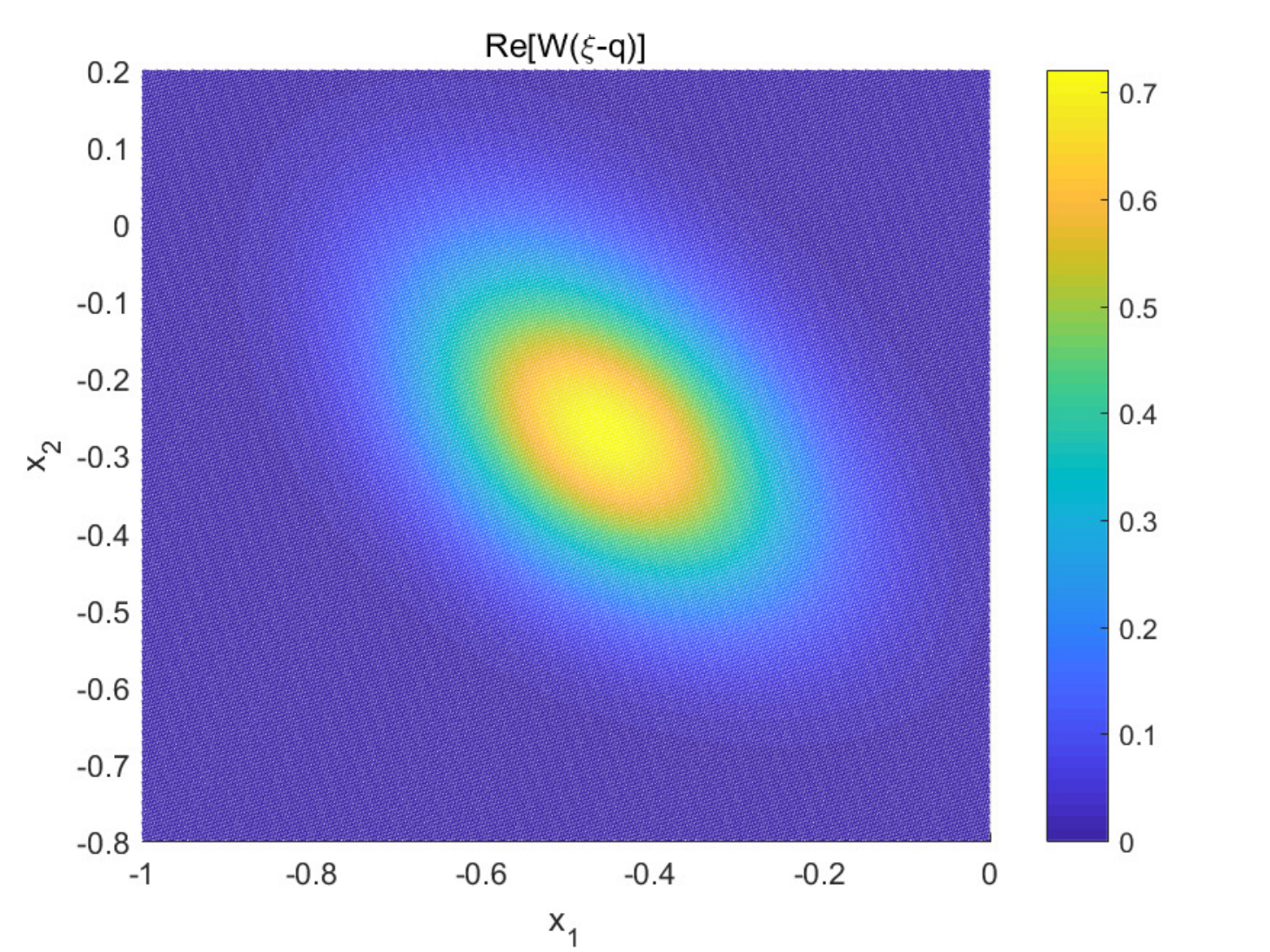} 
\includegraphics[scale=0.45]{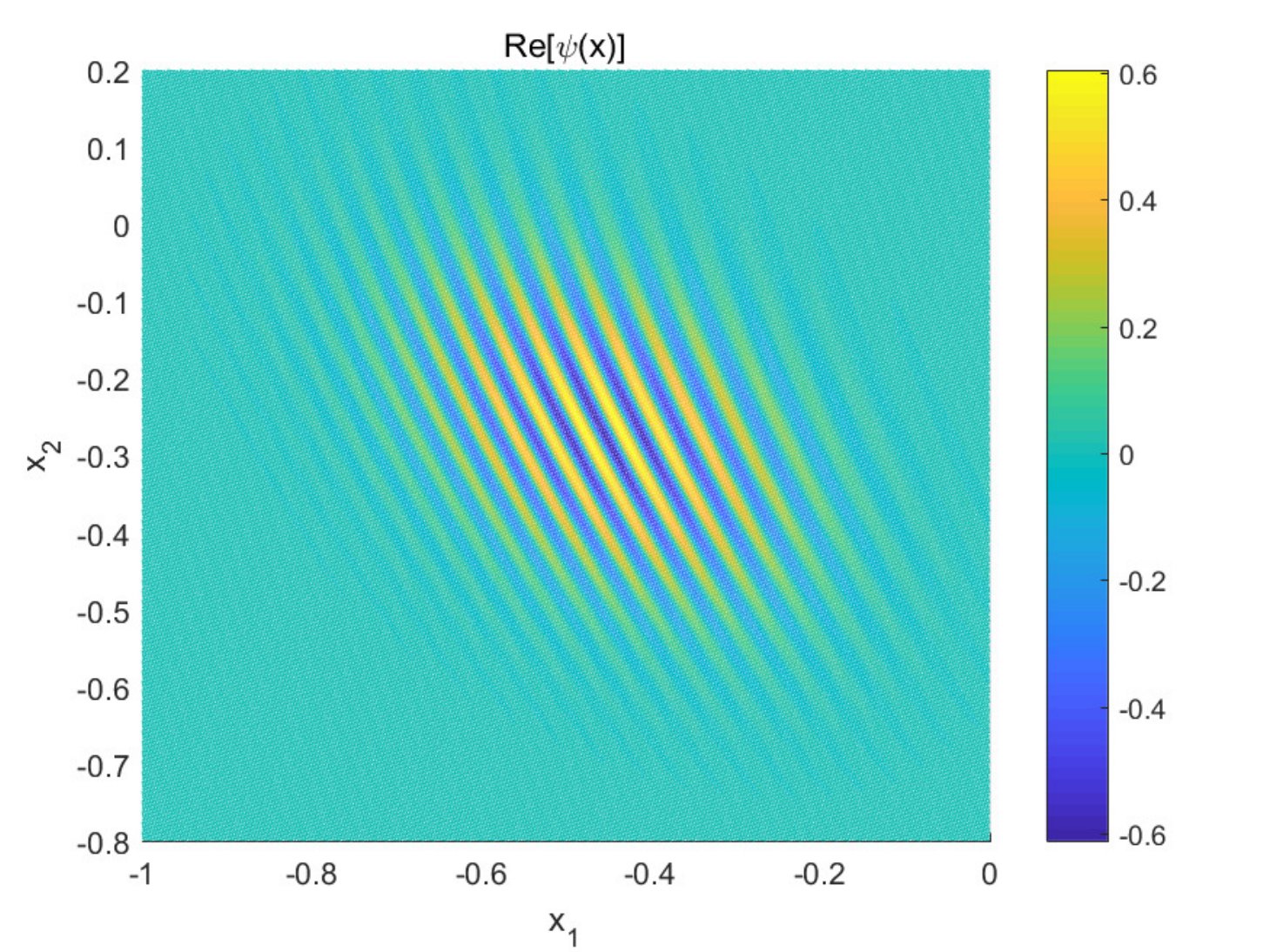} \par
\end{centering}
\caption{
(Example 3)  GWT based SL-TS3 method, $T=0.5$, $\e=\frac{1}{64}$. Left: real part of $W(\bm \xi - \bm q)$. Right: real part of $\psi$.}
\label{fig:2d}
\end{figure} 


Note that, since with resolved spatial mesh, the error in spatial discretization is minimal compared to the error in time. Also, since the time step $\delta t$ for the ODE system is independent of the choice of $\Delta t$, and  the ODE's are cheaper to solve, we choose $\delta t=\Delta t/40$, which makes the error in evolving the parameters negligible. Thus, one expects the numerical error is dominated by the error in time discretization. We plot the results in Table \ref{table:time2d}, which clearly shows the second order convergence in $\Delta t$. 
We remark that, the same tests have been done with the SL-TS2 method, and the results are almost the same, so we omit this part in the paper.

\begin{table}
\small
  \centering
  \begin{tabular}{ c|c | c| c|c} \hline
    Relative error &  $\Delta t= \frac{1}{8}$ & $ \Delta t= \frac{1}{16}$ & $\Delta t= \frac{1}{32}$ & $\Delta t= \frac{1}{64}$    \\ \hline
$\varepsilon= \frac{1}{16}$ & 6.103e-3 & 1.569e-3 & 4.241e-4 & 1.322e-4 
    \\ \hline
 $\varepsilon=\frac{1}{32}$ & 6.145e-3 & 1.624e-3 & 4.707e-4  & 1.616e-4 
    \\ \hline
 $\varepsilon=\frac{1}{64}$ & 6.083e-3 & 1.154e-3 & 4.030e-4  & 1.148e-4 
    \\ \hline
\end{tabular}
\caption{(Example 3) $\varepsilon=\frac{1}{16}$, $\frac{1}{32}$ and $\frac{1}{64}$. Reference solution SL-TS: $\Delta x=\frac{2\pi \eps}{32}$, $\Delta t=\frac{\eps}{32}$. Comparing with GWT based SL-TS3 method: $\Delta \eta=\frac{4\pi}{1024}$,  $\Delta t=\frac{1}{8}$, $\frac{1}{16}$, $\frac{1}{32}$ and $\frac{1}{64}$.}
  \label{table:time2d}
\end{table}

Finally, to explore the dynamical behavior of a quantum wave packet due to the vector potential,  we choose  Gaussian wave packet parameters as
\[
q_0=\left(0.4,0.3\right)^T,\quad p_0=(0,0)^T,\quad C=i \,{\mathrm I}_2,\quad \gamma_0=0,
\]
and we set $\tilde V =0$ and keep the vector potential unchanged. Since the initial wave packet momentum is zero and the scalar potential is zero, the motion of the (classical) wave packet trajectory is solely determined by the initial position and the vector potential. For $\e=\frac 1 8$, $ \frac 1 {16}$ and $\frac 1 {32}$, we compute the expectation of the position with the GWT based SL-TS3 method till $t=2$, which are compared with the classical trajectory in Figure \ref{fig:traj}. Here, the expectation values are computed by \eqref{xave}, namely, they are directly computed using the $w$ function without reconstructing the wave function. We observe clearly from the zoom-in plot that, as $\e \rightarrow 0$, the trajectories of the spatial averages converge to the classical one. We can further examine the snapshots of the wave function $\psi$ and the density function  $\rho=|\psi|.^2$ at $t=0$, $t=1$ and $t=2$ respectively, as plotted in Figure \ref{fig:snap}. First, we observe that initially the wave function does not have an oscillatory phase, but it picks up phase oscillations through dynamics. Second, the wave packet quickly spreads out as it evolves in time, which implies long time simulation based on wave packet approaches might be more challenging.

  \begin{figure}
\begin{centering}
\includegraphics[scale=0.45]{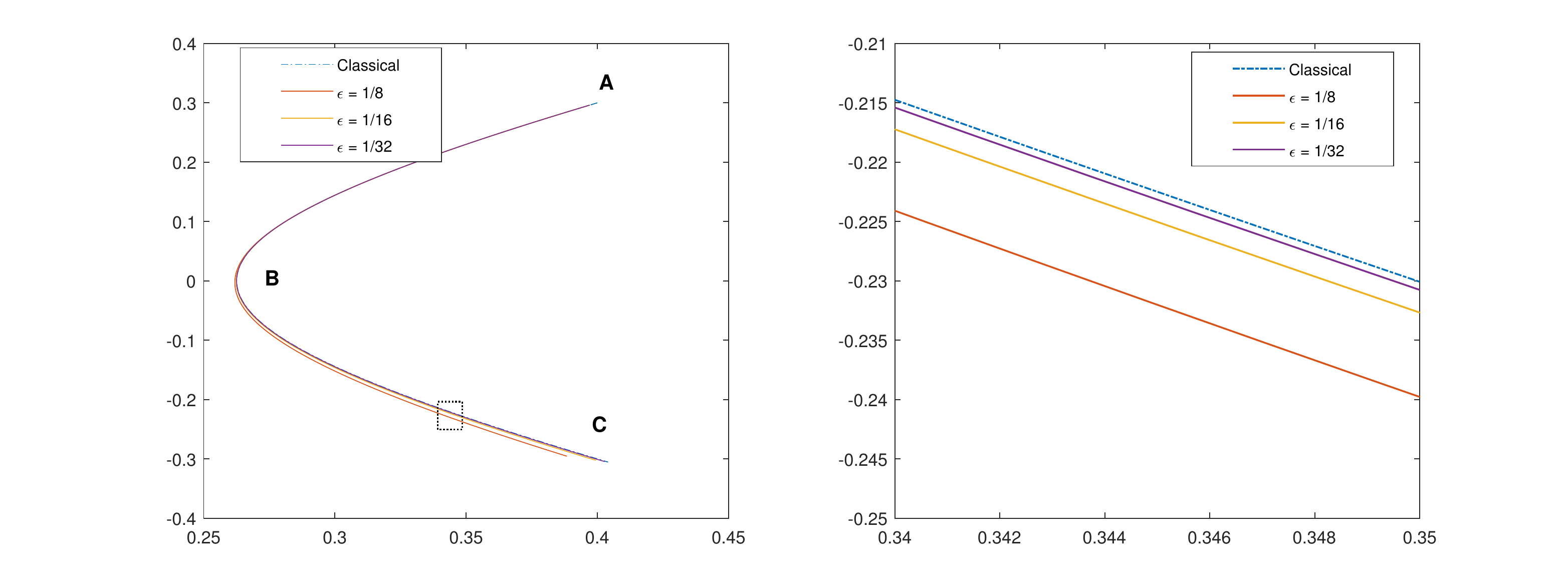}  \par
\end{centering}
\caption{
Classical trajectory (dash-dot line) versus the expectation value of the position (solid lines) for $\e= \frac{1}{8}$, $\frac 1 {16}$ and $\frac 1 {32}$. Left: Overview plot. Right: Zoom-in plot of the dotted rectangle on the left. Three marks on the trajectory: A, $t=0$; B, $t=1$; C, $t=2$.}
\label{fig:traj}
\end{figure} 

\begin{figure}
\begin{centering}
\includegraphics[scale=0.33]{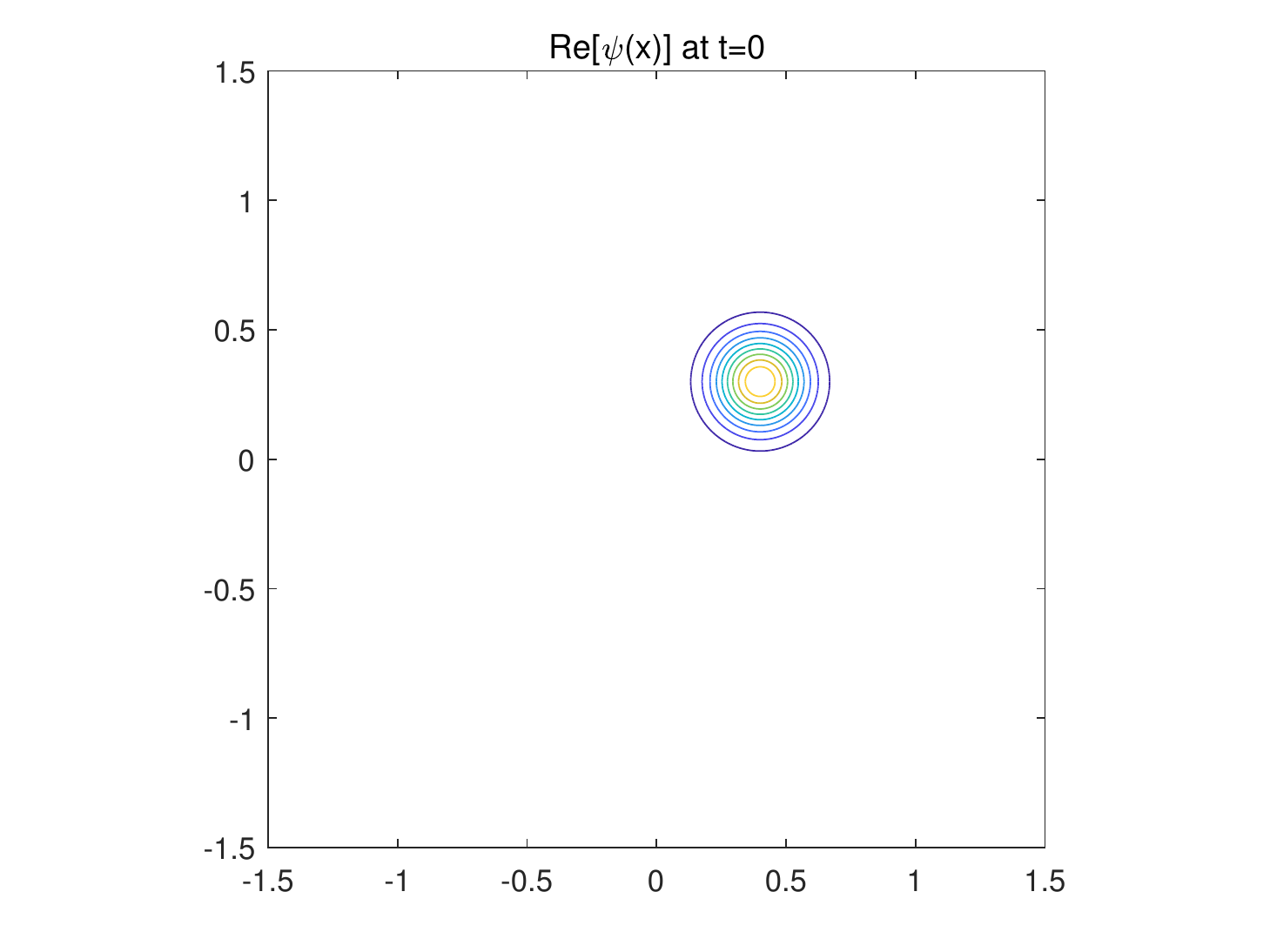}\includegraphics[scale=0.33]{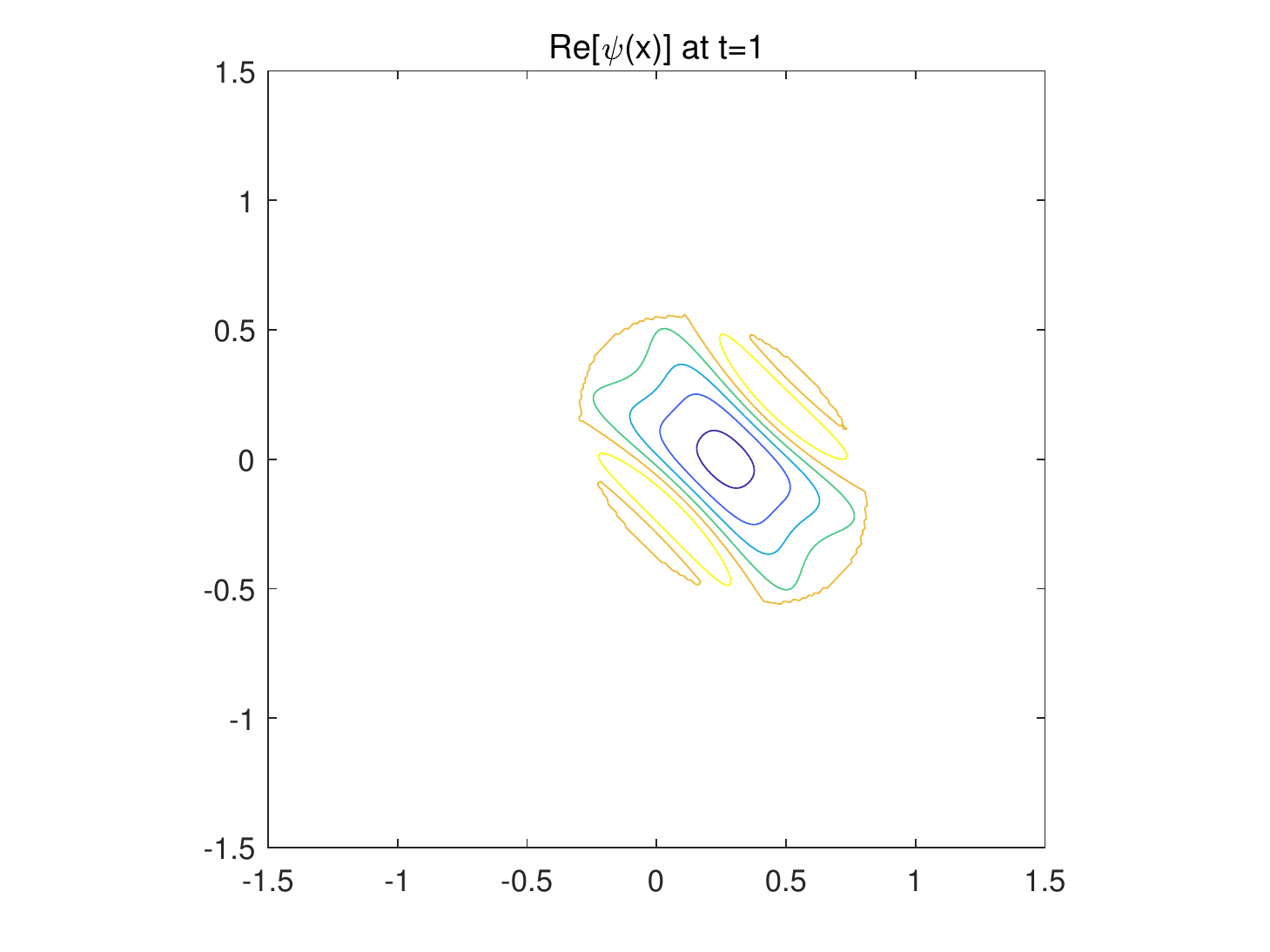}\includegraphics[scale=0.33]{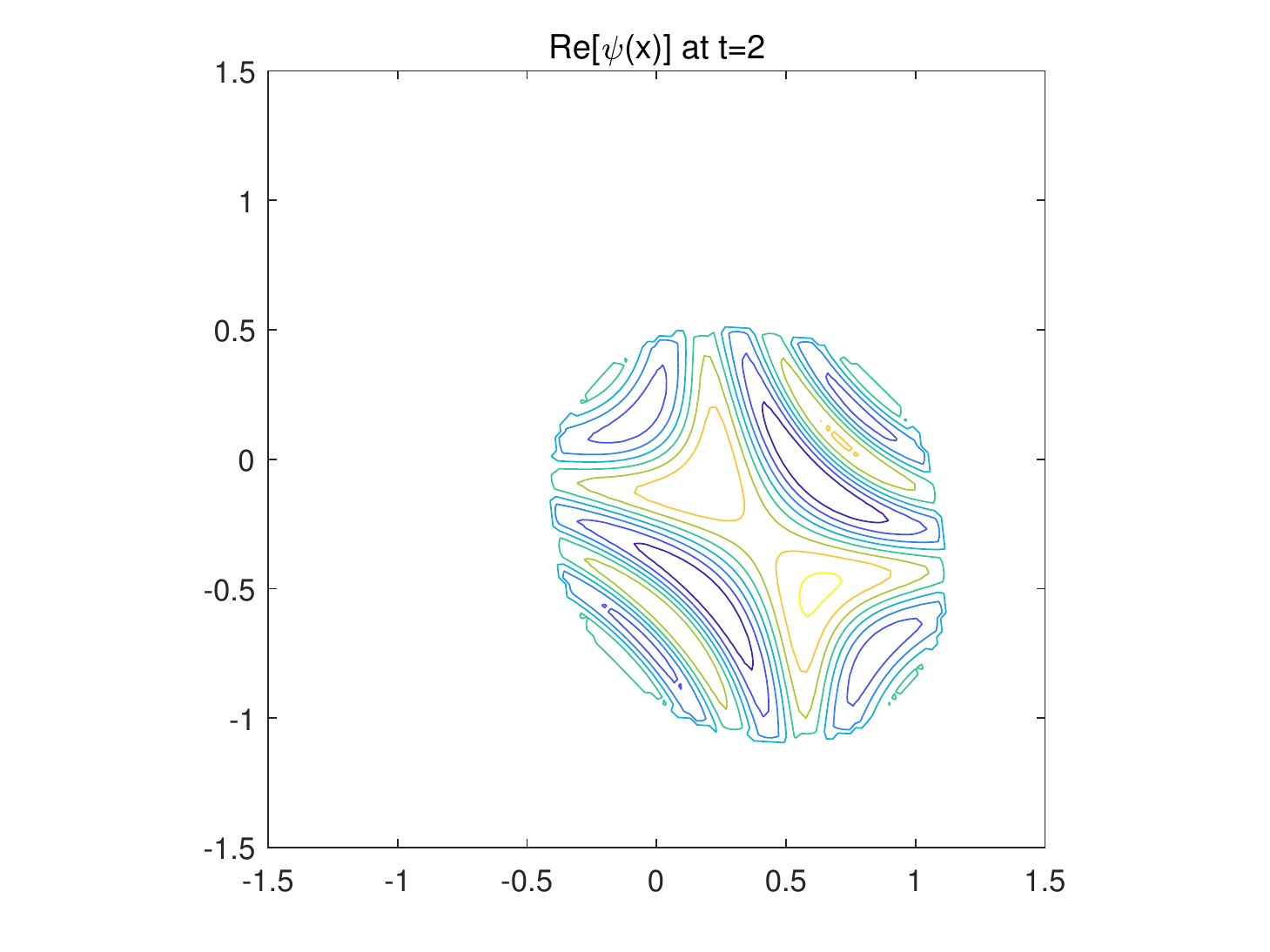}  \par
\includegraphics[scale=0.33]{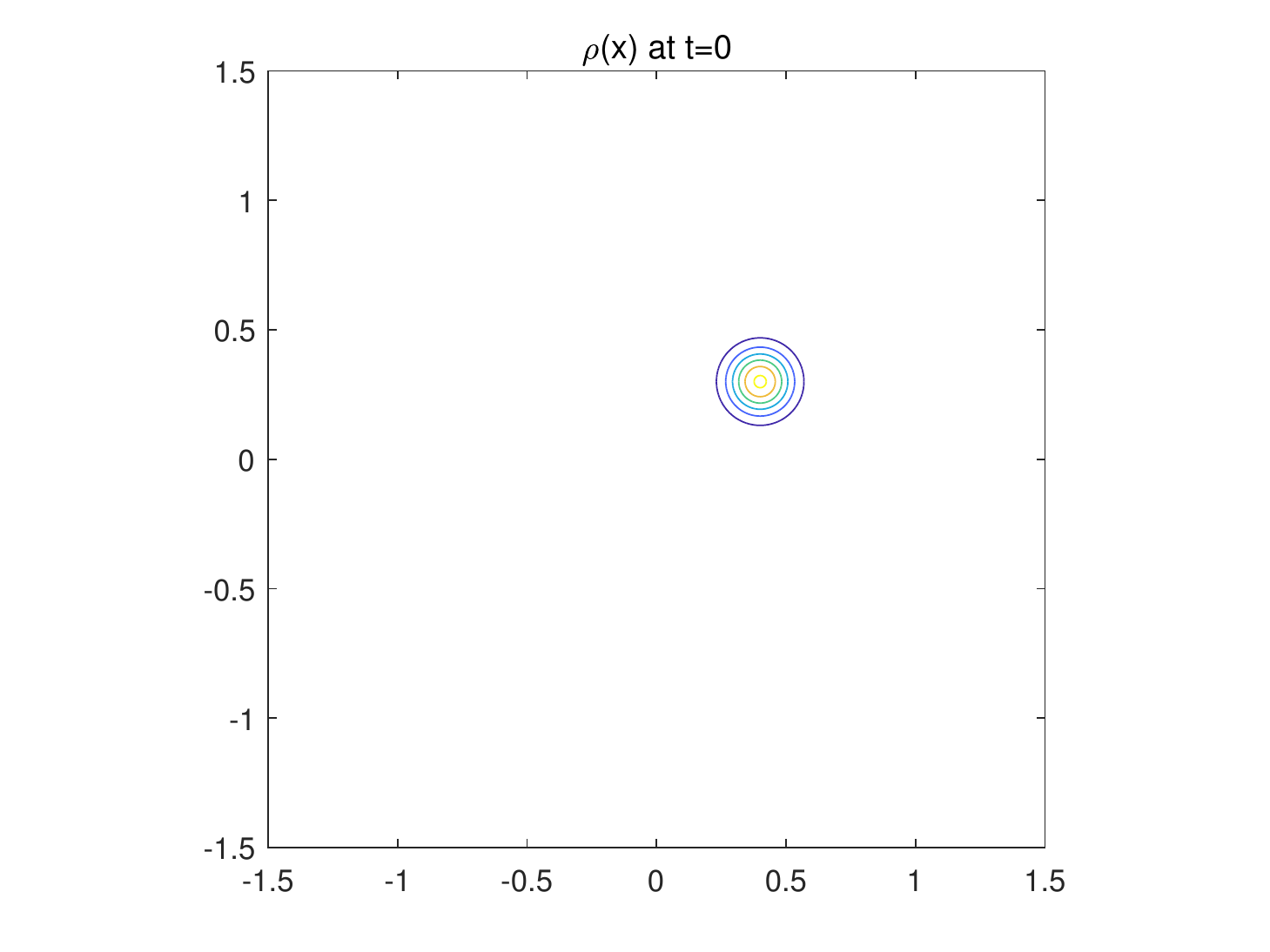}\includegraphics[scale=0.33]{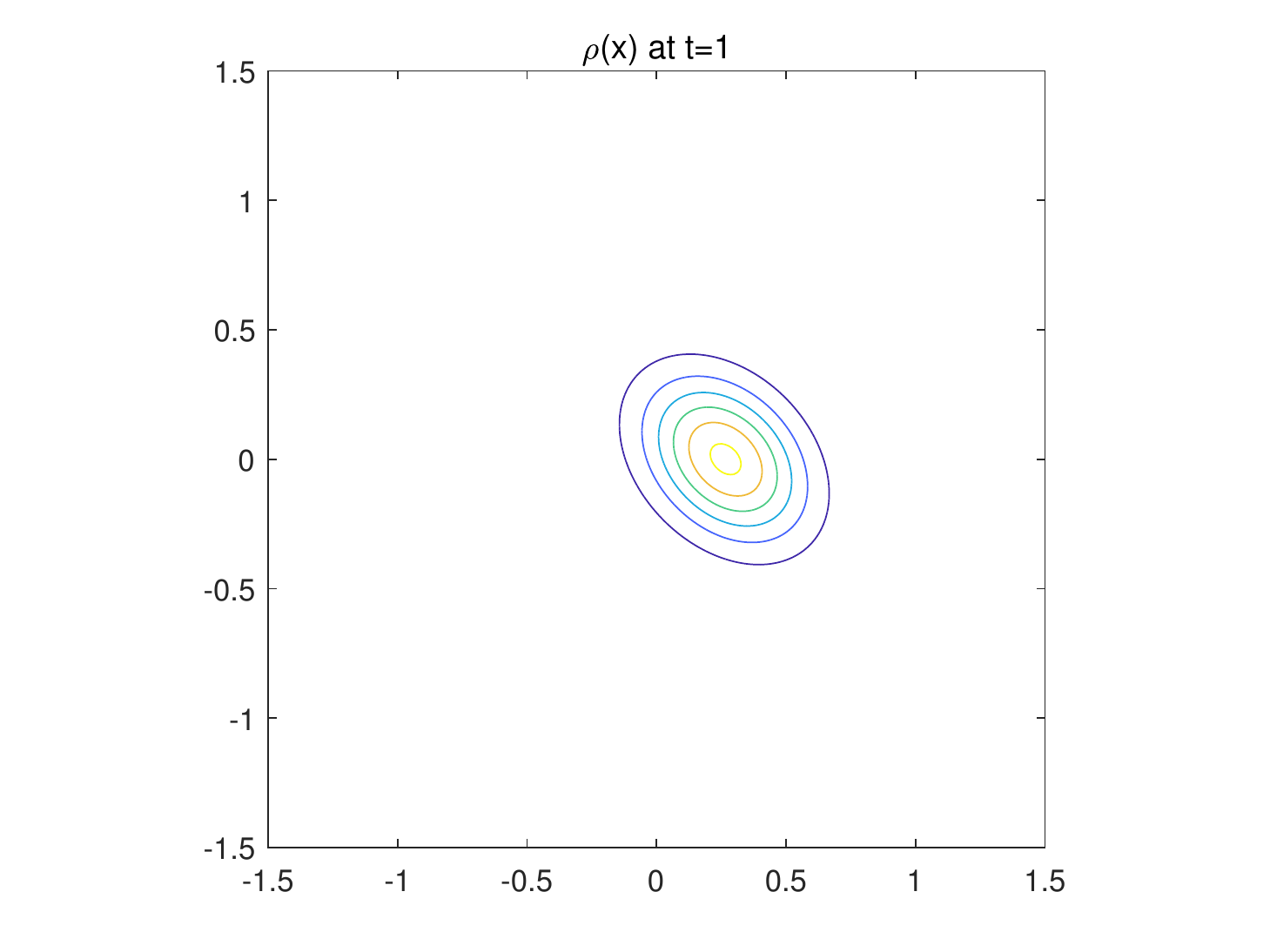}\includegraphics[scale=0.33]{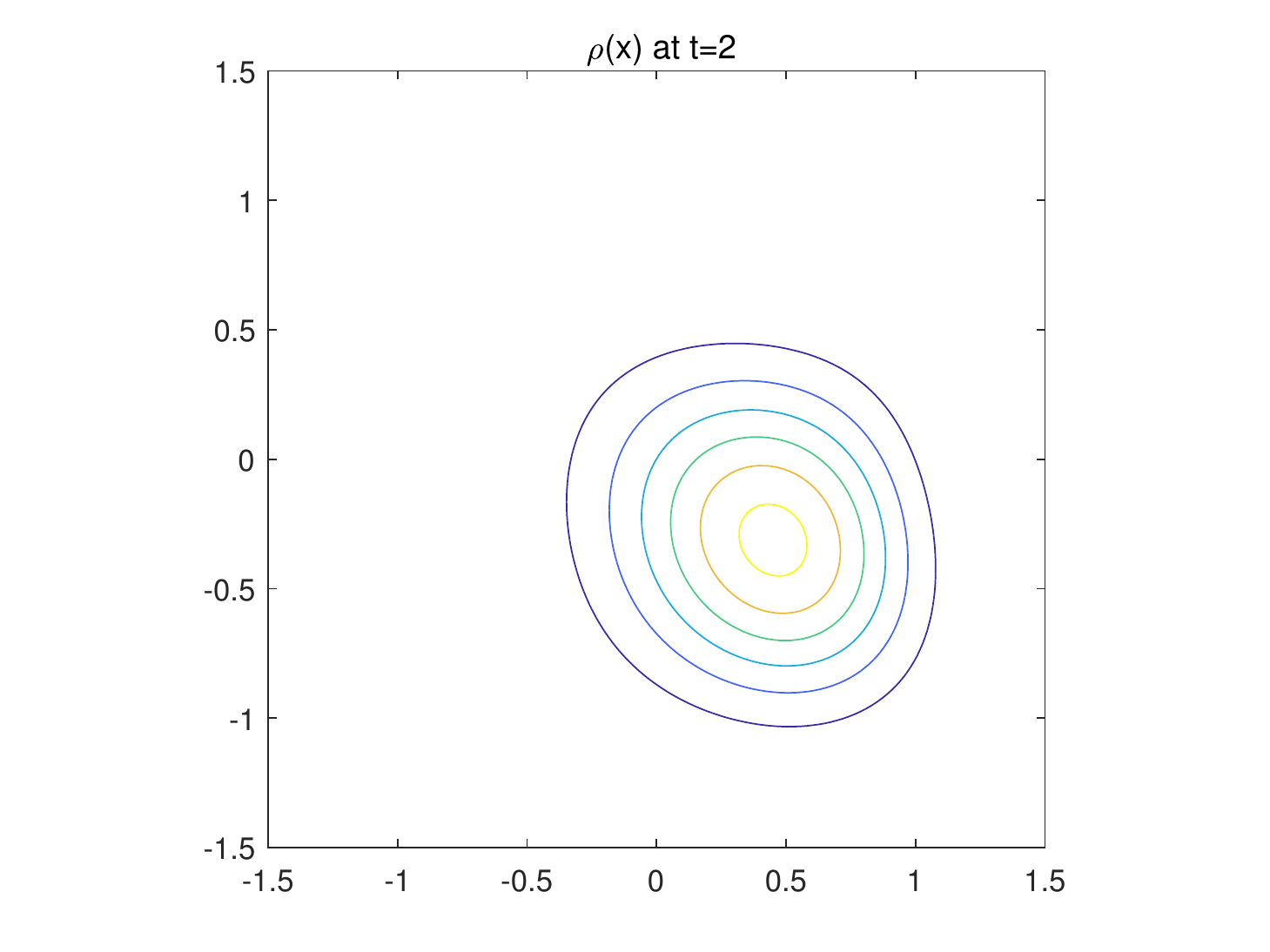} \par
\end{centering}
\caption{
(Example 3) Time snapshots of the real parts of the wave function and density function $\rho=|\psi|^2$ at the three marks in Figure 3. Left: A, $t=0$. Middle: B, $t=1$. Right: C, $t=2$.}
\label{fig:snap}
\end{figure} 

\section{Conclusion}
In this paper we extend the Gaussian wave packet transform method to the Schr\"odinger equation with vector potential, thus being able to describe, among other things, the behavior of a quantum particle in an electromagnetic field, and in particular the small deviation of the expectation value of position and momentum of the particle from the classical counterpart, with great accuracy. 

The new method is a non-trivial generalization of the GWT technique applied to the semi-classical Schr\"odinger equation with scalar potential. 
Indeed, the new approach combines the GWT method with recently developed 
techniques for the treatment of time dependent Schr\"odinger equation with vector potentials (see \cite{SL-TS,NUFFT}), resulting in a very effective strategy for detailed computation of quantum systems in conditions close to the semi-classical limit. 

As in the case of GWT for the Schr\"odinger equation with scalar potential, the method is based on a factorization of the wave function into a fast oscillating phase factor and a slowly varying function, $w$, which satisfies another Schr\"odinger type equation which is more amenable for numerical computation. 
Very few modes are needed to resolve the rescaled wave function $w$, thus allowing accurate direct computation in several space dimensions. Also, the potential in the equation for $w$ is an $O(\sqrt{\eps})$ perturbation of a harmonic oscillator, and therefore, the equation becomes easier to solve when approaching the semi-classical limit. 
Compared to the scalar case, the equation with vector potentials has an additional convection term, which poses some challenge for accurate computation. Such a problem has been successfully dealt with a very effective semi-lagrangian approach, in which Fourier interpolation at the foot or characteristics allows spectral accuracy in space. Such Fourier interpolation requires the computation of the inverse Fourier transform on scattered points. In order to avoid the double summation, we have adopted the Nonuniform FFT algorithms by Greengard (see \cite{NUFFT,nufft6}) to significantly reduce the cost in the spectral interpolation.   

Two splitting methods are presented in the paper: SL-TS2 and SL-TS3, based, respectively, on a two-step and a three-step splitting of the equation. Second order accurate Strang splitting has been used for all the numerical tests reported in the paper. For such a splitting, the two approaches are comparable in term of cost and accuracy, when using the same space resolution and time step. The resulting method is spectrally accurate in space, while only second order actuate in time. 

High order splitting methods can be used, for improving the accuracy. In \cite{GB Trans,GB Trans2},
the fourth order splitting by Chin and Chen in \cite{ChinChen1,ChinChen2} was adopted  with a suitable extension to the Schr\"odinger equation with time dependent kinetic energy. Compared to other more standard fourth order splitting method, the Chin and Chen's method  is particularly accurate and inexpensive, because of the special structure of the Schr\"odinger equation. It would be interesting to investigate the extension of such a method to the case studied in this paper. 

Finally, we observe that in our tests, the error in the numerical solution increases exponentially in time, in agreement with other wave packet based approaches (see, e.g. \cite{Hagedornraising,ErrorestiGB,FGBtrans}).
The phenomenon is however not fully understood and deserves further investigation. In particular, it would be interesting to explore the possibility of following the semi-classical behavior of the system for much longer time.


\section*{Acknowledgments}

The research has been partially funded by 
ITN-ETN Marie-Curie Horizon 2020 program ModCompShock, {\em Modeling and computation of shocks and interfaces}, Project ID: 642768; 
by project F.I.R. 2014 {\em Charge transport in graphene and low dimensional systems}, University of Catania.
Z. Zhou is partially supported by RNMS11-07444 (KI-Net) and the start up grant from Peking University.

\appendix
\appendixpage

\section{Derivation of the Gaussian wave packet transform}
In this section, we present the detailed calculation in the application of the Gaussian wave packet transform to the semi-classical Schr\"odinger equation with vector potentials. Recall that, the Gaussian wave packet transform takes the following ansatz 
\begin{equation}\label{ansats2}
\psi(\mathbf x,t)=W(\bm \xi,t) e^{f(\bm \xi, t)}=W(\bm \xi,t)\exp \left(i\left( \bm \xi^{T} \bm \alpha_R \bm \xi +\bm p^{T}\bm \xi+\gamma_2 \right)/\e \right),
\end{equation}
where $\bm q $, $\bm p$, $\gamma_2$ and $ \bm \alpha_R$ are solely functions of $t$. And for convenience, we recall that $\bm \xi = \bm x -\bm q$ satisfies
\[
\nabla_{\bm x} \bm \xi = \bm I, \quad \partial_t \bm \xi = - \partial_t \bm q, 
\]
where $\bm I$ stands for the identity matrix.

Then, by direct calculations, we have
\begin{align*}
\partial_t \psi & = \partial_t W e^f + W  \partial_t  f e^f \\
& =e^f \left(W_{t}-\nabla_{\bm \xi}W^T \dot{\bm q}\right) + \frac{i}{\e} e^f W\left(\bm \xi^{T}\dot{\bm \alpha_{R}}\bm \xi+\bm \xi^{T}\dot{\bm p}+\dot{\gamma}_{2}-2\bm \xi^{T}\bm \alpha_{R}\dot{\bm q}-\bm p^{T}\dot{\bm q}\right),
\end{align*}
\begin{equation*}
\nabla_{\bm x} \psi  = \nabla_{\bm x} W e^f + W \nabla_{\bm x} f e^f  =e^f \nabla_{ \bm \xi}W +\frac{i}{\varepsilon} e^f \left(2 \bm \alpha_{R} \bm \xi+\bm p\right)W,
\end{equation*}
and
\begin{align*}
\Delta_{\bm x} \psi & =  \Delta_{\bm x} W e^f +2 \nabla_{\bm x} W \cdot \nabla_{\bm x} f e^f + W |\nabla_{\bm x} f|^2 e^f +  W \Delta_{\bm x} f e^f\\
& = \Delta_{\bm \xi} W e^f + \frac{i}{\e} \nabla_{\bm \xi} W^T \left( 4 \bm \alpha_R \bm \xi +2 \bm p \right) e^f +\frac{i}{\e} 2 W {\rm Tr}(\bm \alpha_R) e^f \\
& \quad - \frac{1}{\e^2} W \left( 4 \bm \xi^T \bm \alpha_R^2  \bm \xi +  4 \bm \xi^T \bm \alpha_R \bm p + |\bm p|^2 \right).
\end{align*}


\section{A three dimensional example: magnetic field along the $z$ direction} \label{3dex}

We consider the three dimensional Schr${\rm \ddot o}$dinger equation with a magnetic field $\bm{B} = (0,0,B(x,y))^T$ along the $z$ axis, a harmonic scalar potential in the $z$ direction and an arbitrary scalar potential $V_0(x,y)$, where $B(x,y)$ is the magnified of the magnetic field. We further assume that the scalar potential and the vector potential take the following forms
\[
V= \frac 1 2 z^2+ V_0(x,y),\quad \mathbf{A} =  (A_1 (y), A_2(x), 0)^T.
\]
We remark that, the vector potentials may take other forms but we only consider this specific form for simplicity of calculation.
Obviously, with this choice of the vector potentials we have $\nabla \cdot \mathbf A =0$, and \[
B(x,y)= \frac{d }{dx} A_2 - \frac{d }{d y} A_1 . 
\]
Note that, in this case, the quantum Hamiltonian naturally admits the following decomposition
\[
H=H_{xy}+ H_z,
\]
where
\begin{align*}
H_{xy}&= \frac 1 2 (- i \e \partial_x -A_1 (y) )^2 + \frac 1 2 (- i \e \partial_y - A_2(x))^2+ V_o(x,y) \\
 &= - \frac {\e^2}{2}(\partial_{xx}+ \partial_{yy}) + i \e (A_1 \partial_x  +A_2 \partial_y ) +\tilde V,
\end{align*}
\[
\tilde V= \frac{1}{2}(A_1^2+A_2^2)+ V_o,
\]
and 
\[
H_z=- \frac {\e^2}{2} \partial_{zz} + \frac 1 2 z^2.
\]
Note that $H_z$ corresponds to the well-known quantum harmonic oscillator in the $z$ variable with eigenvalues $E_k=(k+\frac 1 2)\e$ and associated eigenfunctions denoted by $\phi_k(z)$. Due to the completeness and orthogonality of $\phi_k(z)$, if the initial condition is chosen as 
\[
\psi({\bf x},0)=\tilde \psi_0(x,y)\phi_k(z),
\]
then we have 
\[
\psi({\bf x},t)=\tilde \psi (x,y, t)\phi_k(z),
\]
where $\tilde \psi$ satisfies
\begin{equation} \label{eq:2D}
i \e \partial_t \tilde \psi = (H_{xy}+E_k)\tilde \psi.
\end{equation}
Thus, we have reduced the three dimensional quantum dynamics to a two dimensional quantum dynamics and a harmonic oscillator in the $z$ direction. Clearly, if we take $u(x,y,t)=\psi (x,y, t)e^{i E_k t/\e}$, then $u$ satisfies
 \begin{equation} \label{eq:2Du}
i \e \partial_t u = H_{xy} u,  \quad u(x,y,0)=\tilde \psi_0(x,y).
\end{equation}
Hence, for simplicity of simulation, we can  numerically solve \eqref{eq:2Du} for $u(x,y,t)$, and the total wave function is obtained by
\[
\psi({\bf x},t)=\tilde \psi (x,y, t)\phi_k(z)= u(x,y,t)e^{-i E_k t/\e}\phi_k(z).
\]

\end{document}